\input amstex
\documentstyle{amsppt}
\magnification=\magstep1 \NoRunningHeads
\topmatter
\title
Directional recurrence and directional  rigidity  \\ for infinite measure preserving actions \\ of nilpotent lattices
\endtitle

\author
Alexandre I. Danilenko
\endauthor

\thanks
A part of this work was done during a visit of the author to MPIM, Bonn,  in July 2014. The author thanks MPIM for the support and hospitality.
\endthanks

\email
alexandre.danilenko@gmail.com
\endemail

\address
 Institute for Low Temperature Physics
\& Engineering of National Academy of Sciences of Ukraine, 47 Lenin Ave.,
 Kharkov, 61164, UKRAINE
\endaddress
\email alexandre.danilenko\@gmail.com
\endemail

\abstract
Let $\Gamma$ be a lattice in a simply connected nilpotent Lie group $G$.
Given an infinite  measure preserving action $T$ of $\Gamma$ and a ``direction'' in $G$ (i.e. an element $\theta$  of the projective space $P(\goth g)$ of the Lie algebra $\goth g$ of $G$),   some notions of
recurrence and  rigidity  for $T$ along $\theta$ are introduced.
It is shown that the set of recurrent directions $\Cal R(T)$ and the set of rigid directions  for $T$ are  both $G_\delta$.
In the case where $G=\Bbb R^d$ and $\Gamma=\Bbb Z^d$, we prove that (a) for each $G_\delta$-subset $\Delta$ of $P(\goth g)$ and a countable subset $D\subset\Delta$,
there is  a rank-one  action $T$  such that $D\subset\Cal R(T)\subset\Delta$ and (b) $\Cal R(T)=P(\goth g)$ for a generic infinite measure preserving action $T$ of $\Gamma$.
This  answers partly a question from a recent paper by A.~Johnson and A.~{\c S}ahin.
Some applications to the directional entropy of Poisson actions are discussed.
In the case where $G$ is the Heisenberg group $H_3(\Bbb R)$ and $\Gamma=H_3(\Bbb Z)$, a rank-one  $\Gamma$-action $T$ is constructed for which $\Cal R(T)$ is not invariant under the natural ``adjoint'' $G$-action.
\endabstract

\endtopmatter

\document

\head 0. Introduction
\endhead
Subdynamics is the study of the relationship between the dynamical properties of the action of a group
$G$, and those of the action restricted to subgroups of
$G$.
In this paper we consider  measure preserving actions defined on $\sigma$-finite standard measure spaces.
 In the 1980's Milnor generalized
the study of sub-dynamics by defining a concept of {\it directional  entropy} of a
$\Bbb Z^d$-action in every direction,
including the irrational directions for which there is no associated subgroup action \cite{Mi}.
To this end  he considered  $\Bbb Z^d$ as a lattice in $\Bbb R^d$ and he exploited  the geometry of  mutual  position of this lattice and  the 1-dimensional subspaces (i.e. directions) in $\Bbb R^d$.
For a detailed account on the directional entropy of $\Bbb Z^2$-actions and some applications to topological dynamics (expansive subdynamics) we refer to \cite{Pa} and references therein.
In a recent paper \cite{JoSa}, Johnson and {\c S}ahin applied the ``directional approach''
to study {\it  recurrence properties} of infinite measure preserving $\Bbb Z^2$-actions.
They were motivated by 
Feldman's proof of the ratio ergodic theorem \cite{Fel}.
In particular, they showed that for each such an action, say $T$,  the set  $\Cal R(T)$ of all recurrent directions of $T$ is a $G_\delta$-subset of the circle $\Bbb T$.
They also exhibited  examples of   rank-one actions $T$ and $T'$ with  $R(T)=\emptyset$ and $\Bbb T\ne\Cal R(T')\supset \{e^{\pi iq}\mid q\in\Bbb Q\}$. 
They  raised a question: which $G_\delta$-subsets of $\Bbb T$ are realizable as recurrence sets, i.e. appear as $R(T)$ for some $T$? 
We answer this question in part.
\roster
\item"---"
We show that {\it each countable} $G_\delta$ is a recurrence set.
\item"---"
More generally, for each $G_\delta$-subset  $\Delta$ of the projective space $P(\Bbb R^d)$ and a countable subset $D$ of $\Delta$, there is a rank-one infinite measure preserving free $\Bbb Z^d$-action $T$ such that $D\subset\Cal R(T)\subset\Delta$ (Theorem~4.2).
\item"---"
We also prove that  a {\it generic} infinite measure preserving action $T$ of $\Bbb R^d$ is recurrent in every direction, i.e. $\Cal R(T)=P(\Bbb R^d)$ (Theorem~5.6).
\endroster
In parallel to this we introduce a concept of directional rigidity for $\Bbb Z^d$-actions and obtain similar results for realization of $G_\delta$-subsets of 
$P(\Bbb R^d)$ as rigidity sets.

As a byproduct, we obtain some examples of Poisson $\Bbb R^d$-actions with the following  entropy properties:
\roster
\item"---" There is a Poisson action $V=(V_g)_{g\in\Bbb R^d}$ of $0$
 entropy such that for each non-zero $g\in\Bbb R^d$, the transformation $V_g$ is Bernoullian of infinite entropy (Proposition~5.7).
 \item"---"
  For each $G_\delta$-subset  $\Delta\subset P(\Bbb R^d)$ and a countable subset $D$ of $\Delta$, there is a 
   Poisson action $V=(V_g)_{g\in\Bbb R^d}$ of $0$
 entropy such that
 for each nonzero $g\not\in\bigcup_{\theta\in\Delta}\theta$, the transformation
$V_g$ is Bernoulli  of infinite entropy and for each $g\in \bigcup_{\theta\in D}\theta$, the transformation $V_g$ is rigid and hence of 0 entropy (Proposition~5.8).
 \endroster
 In this connection we recall the main result from 
\cite{FeKa}: there is  a Gaussian action $V=(V_g)_{g\in\Bbb Z^2}$  of   $0$ entropy such that every transformation $V_g$, $0\ne g\in\Bbb Z^2$, is Bernoullian.

We extend the concepts of directional recurrence and directional rigidity to actions of lattices $\Gamma$ in simply connected nilpotent  Lie groups $G$.
By a ``direction'' we now mean a 1-parameter subgroup in $G$.
Thus the set of all directions is the projective space $P(\goth g)$, where $\goth g$ denotes the Lie algebra of $G$.
As in the Abelian case (considered originally in \cite{JoSa}), we show that
\roster
\item"---"
Given a measure preserving action $T$ of $\Gamma$, the set  $\Cal R(T)$ of all
recurrent directions of $T$ is a $G_\delta$ in $P(\goth g)$ (Theorems 2.5 and 2.6).
\endroster
Since $G$ acts on $P(\goth g)$ via the adjoint representation, we define another invariant   $\Cal E\Cal R(T)$ of {\it even recurrence} for $T$ as the largest $G$-invariant subset of $\Cal R(T)$.
\roster
\item"---"
Some examples of rank-one actions $T$ of the Heisenberg group $H_3(\Bbb Z)$ are constructed for which $\Cal R(T)$ is either empty (Theorem~6.1) or countably infinite (Theorem~6.2) or uncountable (Theorem~6.3)\footnote{We consider $H_3(\Bbb Z)$ as a lattice in the 3-dimensional real Heisenberg group $H_3(\Bbb R)$.}.
\item"---"
An example of $T$ is given such that $\Cal E\Cal R(T)\ne\Cal R(T)$ (Theorem~6.2).
\endroster
Given an action $T$ of $\Gamma$, we can define a natural analog of the ``suspension flow'' corresponding to $T$.
This is the {\it induced} (in the sense of Mackey) action $\widetilde T$
of $G$ associated with $T$.
Since $\Cal R(T)$ coincides with  the set $\Cal R(\widetilde T)$ of conservative  $\Bbb R$-subactions of $\widetilde T$ in the Abelian case \cite{JoSa}, it is natural to conjecture that 
$\Cal E\Cal R(T)=\Cal R(\widetilde T)$ in the genaral case.
It remains an open problem.
However the analogous claim for the rigidity  sets does not hold in the non-Abelian case~(Remark~2.2).

The outline of the paper is as follows.
In Section~1 we introduce the main concepts and invariants related to the directional recurrence and rigidity.
In Section~2 we discuss relationship between the directional recurrence and rigidity of an action of a lattice in a nilpotent Lie group and similar properties of the {\it suspension flow}, i.e. the induced action of the underlying Lie group.
It is also shown there that the sets of recurrent and rigid directions are both $G_\delta$.
In Section~3 we recall the $(C,F)$-construction of rank-one actions and provide a sufficient condition for  directions to be recurrent in terms of the $(C,F)$-parameters.
This condition is used in Section~4 to construct rank-one actions of $\Bbb Z^d$ with various  sets of recurrent directions.
In Section~5 we prove that a generic $\Bbb Z^d$-action is recurrent in every direction.
This section contains also some applications to the directional entropy of Poisson actions.
In Section~6 we study directional recurrence of infinite measure preserving actions of $H_3(\Bbb Z)$.
The final Section~7 contains a list of open problems and concluding remarks.

{\it Acknowledgements.}  I thank T.~Meyerovitch and  E.~Roy for useful discussions   concerning the  entropy of Poisson suspensions.

\head 1. Recurrence, even recurrence, rigidity and even rigidity along directions
\endhead

Let $G$ be a simply connected   nilpotent Lie group, $\goth g$ the Lie algebra of $G$ and $\exp:\goth g\to G$ the exponential map.
We note that $\exp$ is  a diffeomorphism of $\goth g$
onto $G$ \cite{Mal}.
Let $P(\goth g)$ denote the set of lines (i.e. 1-dimensional subspaces) in $\goth g$.
We endow $P(\goth g)$ with the usual topology of projective space.
Then $P(\goth g)$ is a compact manifold.
The adjoint $G$-action on $\goth g$ induces a smooth $G$-action on $P(\goth g)$.
We denote this action by the symbol ``$\cdot$''.
Given $v\in \goth g\setminus\{0\}$, we let $\exp(v):=\{\exp(tv)\mid {t\in\Bbb R}\}$.
 Then $\exp(v)$ is a closed 1-dimensional subgroup of $G$.
 We note that if $w=tv$ for some $t\in\Bbb R\setminus\{0\}$ then $\exp(w)=\exp(v)$.
Hence for each line $\theta\in P(\goth g)$, the notation $\exp(\theta)$ is well defined.
Moreover, $g\exp(\theta)g^{-1}=\exp(g\cdot\theta)$ for each $g\in G$.

Let $R=(R_g)_{g\in G}$ be a measure preserving action of $G$ on a $\sigma$-finite standard measure space $(Y,\goth Y,\nu)$.

\definition{Definition 1.1}
\roster
\item"(i)"
We recall that $R$ is called {\it conservative} if for each  subset $B\in\goth Y$, $\nu(B)>0$, and a compact $K\subset G$,
 there is  an element $g\in G\setminus K$,  such that
 $$
 \nu(B\cap R_gB)>0.
 $$
 \item"(ii)"
 We call $R$ {\it recurrent along a  line $\theta\in P(\goth g)$} if the flow $(\exp(tv))_{t\in\Bbb R} $ is conservative for some (and hence for each) $v\in\theta\setminus\{0\}$.
 \item"(iii)"
We recall that $R$ is called {\it rigid} if there is a sequence $(g_n)_{n\ge 1}$ of elements in $G$ such that $g_n\to\infty$ and
 $$
 \lim_{n\to\infty}\nu(B\cap R_{g_n}B)=\mu(B)
 $$
 for each subset $B\in\goth Y$ of finite measure.
  \item"(iv)"
  We call $R$
  {\it rigid along a  line $\theta\in P(\goth g)$} if the flow $(\exp(tv))_{t\in\Bbb R} $ is rigid for some (and hence for each) $v\in\theta\setminus\{0\}$.

 \endroster
 \enddefinition

  Denote by $\Cal R(R)$ the set of all $\theta\in P(\goth g)$ such that  $R$ is recurrent along $\theta$.
  Denote by $\Cal Ri(R)$ the set of all $\theta\in P(\goth g)$ such that $R$ is rigid along $\theta$.
  Of course, $\Cal Ri(R)\subset\Cal R(R)$.
  It is easy to see that if a $G$-action $R'$ is  isomorphic to $R$ then $\Cal R(R')=\Cal R(R)$
  and $\Cal Ri(R')=\Cal Ri(R)$.

\proclaim{Proposition 1.2}
The sets $\Cal R(R)$  and $\Cal Ri(R)$ are  $G$-invariant.
 \endproclaim
 \demo{Proof}
 Let $\theta\in\Cal R(R)$.
 Fix
 an element $g_0\in G$.
 Take a subset $B\subset Y$ of positive measure and a compact $K\subset G$.
  Since $R$ is recurrent along $\theta$,  there is  $g\in \exp(\theta)$ such that $g\not\in K$ such that $\nu(B\cap R_gB)>0$.
 Hence
 $$
 0<\nu(R_{g_0}B\cap R_{g_0}R_gB)=\nu(R_{g_0}B\cap R_{g_0gg_0^{-1}}R_{g_0}B).
 $$
 Since $g_0gg_0^{-1}\in \exp(g_0\cdot \theta)$ and $g_0gg_0^{-1}\not\in g_0Kg_0^{-1}$, it follows that the flow $(R_g)_{g\in g_0\cdot\theta}$ is conservative.
 Thus $\Cal R(R)$ is $G$-invariant.
 In a similar way we can verify that $\Cal Ri (R)$ is $G$-invariant.
 \qed
 \enddemo

From now on we fix a lattice $\Gamma$  in $G$.
We recall that there exists a lattice in $G$ if and only if the structural constants of $\goth g$ are all rational \cite{Mal}.
Moreover, every lattice in $G$ is uniform \cite{Mal}, i.e. co-compact.
We fix   a right-invariant metric  dist$(.,.)$ on $G$ compatible with the topology.

Let  $T=(T_\gamma)_{\gamma\in\Gamma}$ be a  measure preserving action of $\Gamma$ on  a  $\sigma$-finite  standard measure space $(X,\goth B,\mu)$.

\definition{Definition 1.3}
\roster
\item"(i)"
We call   $T$  {\it recurrent   along a line}  $\theta\in P(\goth g)$  if for each $\epsilon>0$ and every subset $A\in\goth B$, $\mu(A)>0$, there are an element $\gamma\in\Gamma\setminus\{1_\Gamma\}$ and an element $g\in\exp(\theta)$ such that   $\text{dist}(\gamma,g)<\epsilon$ and  $\mu(A\cap T_\gamma A)>0$.
\item"(ii)"
We call   $T$  {\it evenly recurrent   along a line}  $\theta\in P(\goth g)$
 if  $T$ is recurrent along every line from  the $G$-orbit of $\theta$.
 \item"(iii)"
We call   $T$  {\it rigid   along a line}  $\theta\in P(\goth g)$  if
there is a sequence $(\gamma_n)_{n\ge 1}$ of elements in $\Gamma$ such that
$\lim_{n\to\infty}\inf_{g\in\exp(\theta)}\text{dist}(\gamma_n,g)= 0$ and
  $$
  \lim_{n\to\infty}\mu(A\cap T_{\gamma_n} A)=\mu(A)
  $$
  for each subset
 $A\in\goth B$ with $\mu(A)<\infty$\footnote{This means that $T_{\gamma_n}\to\text{Id}$ as $n\to\infty$ in the weak topology on the group of all $\mu$-reserving invertible transformations of $X$.}.
\item"(iv)"
We call   $T$  {\it evenly rigid   along a line}  $\theta\in P(\goth g)$
 if  $T$ is rigid along every line from  the $G$-orbit of $\theta$.
\endroster
\enddefinition

We denote by $\Cal R(T)$ the set of all
$\theta\in P(\goth g)$ such that  $T$ is recurrent along $\theta$.
We denote by $\Cal Ri(T)$ the set of all
$\theta\in P(\goth g)$ such that  $T$ is rigid along $\theta$.
In a similar way, we denote by $\Cal E\Cal R(T)$ and
$\Cal E\Cal Ri(T)$
the set of all $\theta\in P(\goth g)$ such that $T$ is evenly recurrent  along them and evenly rigid  along them respectively.

Of course, $\Cal R(T)\supset\Cal E\Cal R(T)$, $\Cal Ri(T)\supset\Cal E\Cal Ri(T)$,
$\Cal R(T)\supset\Cal Ri(T)$  and  $\Cal E\Cal R(T)\supset\Cal E\Cal Ri(T)$.
For $G$ Abelian,  $\Cal R(T)=\Cal E\Cal R(T)$ and
$\Cal Ri(T)=\Cal E\Cal Ri(T)$.
However, in general $\Cal R(T)\ne\Cal E\Cal R(T)$ (see Theorem~6.2 below)
and $\Cal Ri(T)\ne\Cal E\Cal Ri(T)$.

\remark{Remark \rom{1.4}}
\roster
\item"(i)"
It is easy to see that if $\theta$ is ``rational'', i.e. the intersection $\Gamma\cap\exp(\theta)$ is nontrivial, say there is $\gamma_0\ne 1_\Gamma$ such that $\Gamma\cap\exp(\theta)=\{\gamma_0^n\mid n\in\Bbb Z\}$, then $\theta$ is recurrent if and only if $\gamma_0$ (i.e. the action of $\Bbb Z$ generated by $\gamma_0$) is conservative.
In a similar way, if $\theta$ is rigid if and only if  $\gamma_0$ is rigid.
\item"(ii)"
If  $\theta\in\Cal R(T)$ then we have $\{\gamma\cdot\theta\mid \gamma\in\Gamma\}\subset\Cal R(T)$.
In a similar way, if  $\theta\in\Cal Ri(T)$ then we have $\{\gamma\cdot\theta\mid \gamma\in\Gamma\}\subset\Cal Ri(T)$.
This can be shown in a similar way as in Proposition~1.2 (plus the fact that diet is right-invariant).
\endroster
\endremark

Given $g\in G$ and $\theta\in P(\goth g)$, we denote by $\text{dist}(g,\exp(\theta))$ the distance from $g$ to the closed subgroup $\exp(\theta)$, i.e. 
$$
\text{dist}(g,\exp(\theta)):=\inf_{h\in\exp(\theta)}\text{dist}(g,h)=\min_{h\in\exp(\theta)}\text{dist}(g,h).
$$

Since in  Definition~1.3(i), there is no any estimation (from below) for the ratio $\mu(A\cap T_\gamma A)/\mu(A)$, the following lemma---which is equivalent to Definition~1.3(i)---is more useful for applications.

\proclaim{Lemma 1.5}
Let $\theta\in\Cal R(T)$.
Then for each $\epsilon>0$, a compact $K\subset G$ and a subset $A\subset X$ of finite measure, there is a  Borel  subset $A_0\subset A$ and Borel one-to-one map  $R:A_0\to A$ and a Borel map $\vartheta:A_0\ni x\mapsto\vartheta_x\in\Gamma\setminus K$ such that $\mu(A_0)\ge 0.5\mu(A)$ and $Rx=T_{\vartheta_x}x$ and $\text{\rom{dist}}(\vartheta_x,\exp(\theta))<\epsilon$ for all $x\in A_0$.
\endproclaim
\demo{Proof}
We use a standard exhaustion argument.
Let
$$
\Gamma_\epsilon:=\{\gamma\in\Gamma\setminus\{1\}\mid \text{dist}(\gamma,\exp(\theta))<\epsilon\}.
$$
Enumerate the elements of $\Gamma_\epsilon$, i.e. let $\Gamma_\epsilon=\{\gamma_n\}_{n\ge 1}$.
We now  set $A_1:=A\cap T_{\gamma_1}^{-1}A$, $B_1:=T_{\gamma_1}A_1$, $A_2:=(A\setminus(A_1\cup B_1))\cap  T_{\gamma_2}^{-1}(A \setminus(A_1\cup B_1))$,  $B_2:= T_{\gamma_2}A_2$ and so on.
Then we obtain two sequences $(A_n)_{n\ge 1}$ and $(B_n)_{n\ge 1}$ of Borel subsets of $A$
such that $A_i\cap A_j=B_i\cap B_j=\emptyset$ whenever $i\ne j$ and $T_{\gamma_i}A_i=B_i$ for all $i$.
We let $A_0:=\bigsqcup_{i\ge 1}A_i$ and $B_0:=\bigsqcup_{i\ge 1}B_i$.
It follows from  Definition~1.3(i) that $\mu(A\setminus (A_0\cup B_0))=0$.
Since $\mu(A_0)=\mu (B_0)$, it follows that  $\mu(A_0)\ge 0.5\mu(A)$.
It remains to let $\vartheta_x:=\gamma_i$ for all  $x\in A_i$, $i\ge 1$.
\qed
\enddemo

\head 2. Recurrence and  rigidity~along~directions in terms of the induced $G$-actions
\endhead

Denote by $\widetilde T=(\widetilde T_g)_{g\in G}$ the action of $G$ {\it induced from $T$} (see \cite{Ma}, \cite{Zi}).
We recall that the space of $\widetilde T$ is the product space $(G/\Gamma\times X,\lambda\times\mu)$, where $\lambda$ is the unique $G$-invariant probability measure on the homogeneous space $G/\Gamma$.
To define $\widetilde T$ we first choose  a Borel cross-section $s:G/\Gamma\to G$ of the natural projection $G\to G/\Gamma$.
Moreover, we may assume without loss of generality that   $s(\Gamma)=1_G$
and  $s$ is a homeomorphism when restricted to an open neighborhood   of $\Gamma$, this neighborhood is of full measure and the measure of the boundary of the neighborhood is $0$.
Define a Borel map $h_s:G\times G/\Gamma\to\Gamma$ by setting
$$
h_s(g, g_1\Gamma)=s(gg_1\Gamma)^{-1}gs(g_1\Gamma).
$$
Then $h_s$ satisfies the 1-cocycle identity, i.e.
$h_s(g_2,g_1g\Gamma)h_s(g_1,g\Gamma)=h_s(g_2g_1,g\Gamma)$ for all $g_1,g_2,g\in\Gamma$.
We now set for $g,g_1\in G$ and $x\in X$,
$$
\widetilde T_g(g_1\Gamma,x):=(gg_1\Gamma, T_{h_s(g,g_1\Gamma)}x).
$$
Then $(\widetilde T_g)_{g\in G}$ is a measure preserving  action of $G$ on $(G/\Gamma\times X,\lambda\times\mu)$.
We note that the isomorphism class of $\widetilde T$ does not depend on the choice  of $s$.

\proclaim{Theorem 2.1}
Let $G=\Bbb R^d$ and $\Gamma=\Bbb Z^d$, $d\ge 1$.
Then
 $\Cal R(\widetilde T)=\Cal R(T)$ and $\Cal Ri(\widetilde T)=\Cal Ri(T)$.
\endproclaim
\demo{Proof}
We consider the quotient space $G/\Gamma$ as $[0,1)^d$.
Given $g=(g_1,\dots,g_d)\in\Bbb R^d$, we let $[g]=(E(g_1),\dots,E(g_d))$ and $\{g\}:=(F(g_1),\dots, F(g_d))$, where $E(.)$ and $F(.)$ denote the integer part and the fractional part of a real.
If the cross-section $s:[0,1)^d\to \Bbb R^d$ is given by the formula $s(y):=y$ then we have
$h_s(g,y)=[g+y]$ for all $g\in G$ and $y\in [0,1)^d$.

 (A) We first show that $\Cal Ri(T)=\Cal Ri(\widetilde T)$.
 Let $\theta\in\Cal Ri(T)$.
  Then there are $\gamma_n\in\Gamma$ and $t_n\in\theta$ such that $\text{dist}(\gamma_n,t_n)\to 0$ and $T_{\gamma_n}\to\text{Id}_X$ weakly as $n\to\infty$.
  We claim that $\widetilde T_{t_n}\to\text{Id}_{(G/\Gamma)\times X}$ weakly as $n\to\infty$.
Indeed,
let $\epsilon_n:=t_n-\gamma_n$.
Then
$$
\widetilde T_{t_n}(y,x)=(\{t_n+y\},T_{[t_n+y]}x)=(\{\epsilon_n+y\}, T_{\gamma_n}T_{[\epsilon_n+y]}x).
\tag2-1
$$
Since the Lebesgue measure of the subset $Y_n:=\{y\in [0,1)^d\mid \epsilon_n+y\in[0,1)^d\}$ goes to 1 as $n\to\infty$
and $\{\epsilon_n+y\}=y$ and $[\epsilon_n+y]=0$ for all $y\in Y_n$, it follows that $\widetilde T_{t_n}\to\text{Id}_{(G/\Gamma)\times X}$ as $n\to\infty$.
Thus we obtain that $\theta\in\Cal Ri(\widetilde T)$.

Conversely, let $\theta\in\Cal Ri(\widetilde T)$.
Then there are $t_n\in\theta$, $n\in\Bbb N$, such that
$$
\widetilde T_{t_n}\to\text{Id}_{(G/\Gamma)\times X}\text{ weakly as $n\to\infty$.}\tag2-2
$$
It follows from \thetag{2-1} that the sequence of transformations $y\mapsto \{t_n+y\}$ of $G/\Gamma$ converge
to Id$_{G/\Gamma}$ as $n\to\infty$.
This, in turn, implies that there is a sequence $(\gamma_n)_{n\in\Bbb N}$ of elements of $\Gamma$ such that
$\lim_{n\to\infty }\text{dist}(t_n,\gamma_n)= 0$.
Therefore, Lebesgue measure of the subset $\{y\in G/\Gamma\mid [t_n+y]=\gamma_n\}$ converges to 1 as $n\to\infty$.
Now~\thetag{2-1} and \thetag{2-2} yield that $T_{\gamma_n}\to \text{Id}_X$.
Hence $\theta\in\Cal Ri(T)$.

(B) We now  show that \footnote{Though this fact was originally stated in \cite{JoSa}  we give here an alternative proof because, on our opinion, the proof of the inclusion $\Cal R(T)\subset \Cal R(\widetilde T)$ was not completed there.} $\Cal R(T)=\Cal R(\widetilde T)$.
Take $\theta\in\Cal R(T)$.
 Given a subset $A\subset G/\Gamma\times X$ of positive  measure, a compact $K\subset G$ and $\epsilon>0$, we find  two subsets $B\subset X$ and $C\subset G/\Gamma$ of finite positive measure
 such that
 $$
 (\text{Leb}\times\mu)(A\cap (B\times C))>0.99\text{Leb}(B)\mu(C).\tag2-3
 $$
 For $t\in G$, we set $B_t:=\{y\in B\mid t+y\in B\text{ and }[t+y]=0\}$.
Then we find $\epsilon_1>0$ so small that
 $
 \text{Leb}(B_t)>0.5\text{Leb}(B)
 $
   for each $t\in G$ such that $\text{dist}(t,0)<\epsilon_1$.
By Lemma~1.5,
there are elements $\gamma_1,\dots,\gamma_l\in\Gamma$, $t_1,\dots, t_l\in\theta\setminus K$ and pairwise disjoint subsets $C_1,\dots, C_l$ of $C$
such that $\max_{1\le j\le l}\text{dist}(\gamma_j,t_j)<\min(\epsilon,\epsilon_1)$, the sets $T_{\gamma_1}A_1,\dots, T_{\gamma_l}C_l$ are mutually disjoint subsets of $C$ and  $\mu(\bigsqcup_{j=1}^l C_j)>0.4\mu(C)$.
We now let $A':=\bigsqcup_{j=1}^lB_{t_j}\times C_j$.
Of course, $A'$ is a subset of $B\times C$.
We have
$$
\widetilde T_{t_j}(b,c)=(\{t_j+b\}, T_{\gamma_j}c)\subset B\times C\qquad\text{if }b\in B_j,\text{ and } c\in C_j
$$
for each $j=1,\dots,l$.
Moreover, the sets $\widetilde T_{t_j}(B_j\times C_j)$, $j=1,\dots,l$, are pairwise disjoint and
$(\text{Leb}\times\mu)(\bigsqcup_{j=1}^l(B_j\times C_j))>0.2(\text{Leb}\times\mu)(B\times C)$.
It now follows from~\thetag{2-3} that there is $j\in\{1,\dots,l\}$ such that $(\text{Leb}\times\mu)(\widetilde T_{t_j}(A\cap(B_j\times C_j)\cap A)>0$.
Hence $\theta\in\Cal R(\widetilde T)$.

Conversely, let $\theta\in\Cal R(\widetilde T)$.
Given $\epsilon>0$, let $Y=[1/2,1/2+\epsilon)\subset G/\Gamma$.
It is easy to see that
if $gY\cap Y\ne \emptyset$ for some $g\in G$ then $\text{dist}(g,\Gamma)<\epsilon$
and the map $Y\ni y\mapsto[g+y]\in\Bbb Z^d$ is constant.
Let $A$ be a subset of $X$ of finite positive measure.
Then there is $g\in\theta$ such that $\text{dist}(g,0)>100$ and
$$
0<(\text{Leb}\times\mu)((Y\times A)\cap \widetilde T_g (Y\times A))=\text{Leb}(gY\cap Y)\mu(A\cap T_{\gamma}A),
$$
where $\gamma:=[g+y]\in\Gamma$ for all $y\in Y$.
It follows that $\text{dist}(\gamma,\theta)<\epsilon$ and $\gamma\ne 0$.
Hence $\theta\in\Cal R(T)$.
\qed

\remark{\rom{Remark 2.2}}
We note that the equality  $\Cal Ri(\widetilde T)=\Cal E\Cal Ri(T)$ does not hold for non-Abelian nilpotent groups.
Consider, for instance, the case where   $G=H_3(\Bbb R)$ and $H=H_3(\Bbb Z)$ (see Section~6 for their definition).
Let $T$ be  an ergodic action of $H_3(\Bbb Z)$.
We claim that
 $\widetilde T$ is not rigid and hence $\Cal Ri(\widetilde T)=\emptyset$.
 Indeed, if $\widetilde T$ were rigid then the quotient $G$-action by translations on $G/\Gamma$ is also rigid.
 However the latter action is mixing relative to the subspace generated by all eigenfunctions \cite{Au--Ha}.
 On the other hand, there are examples of weakly mixing  $H_3(\Bbb Z)$-actions $T$ such that $\Cal Ri(T)$ contains the line passing through the center \cite{Da3}.
\endremark

\enddemo

\proclaim{Corollary 2.3}
Let $G=\Bbb R^d$ and $\Gamma=\Bbb Z^d$, $d\ge 1$.
 If  an action $T$  of $\Gamma$ is ergodic and extends to an action $\widehat T$ of $G$  on the same measure space where $T$ is defined then $\Cal R(T)=\Cal R(\widehat T)$.
 \endproclaim
\demo{Proof}
It follows from the condition of the corollary that the induced $G$-action $\widetilde T$ is isomorphic to the product $\widehat T\times D$, where $D$ is the natural $G$-action by translations on  $G/\Gamma$ \cite{Zi, Proposition~2.10}.
Since $D$ is finite measure preserving,  $\Cal R(\widehat T\times D)=\Cal R(\widehat T)$ (see Lemma~2.4(ii) below).
It remains to apply Theorem~2.1.
\qed
\enddemo

We leave the proof of the following non-difficult statement to the reader as an exercise.

\proclaim{Lemma 2.4}
 Let  $F=(F_t)_{t\in\Bbb R}$ be a $\sigma$-finite  measure preserving flow and let $S=(S_t)_{t\in\Bbb R}$ be a probability preserving flow.
\roster
\item"\rom(i)"
$F$ is conservative if and only if the transformation $F_1$ is conservative.
\item"\rom(ii)"
$F$ is conservative if and only if the product flow  $(F_t\times S_t)_{t\in\Bbb R}$ is conservative\footnote{A similar claim for transformations (i.e. $\Bbb Z$-actions) is proved in \cite{Aa}.  We note that (ii) follows from that claim and  (i).}.
\item"\rom(iii)"
$F$ is rigid if and only if  $F_1$ is rigid.
\endroster
\endproclaim

We now describe the ``topological type'' of  $\Cal R(T)$ and $\Cal E\Cal R(T)$ as subspaces of $P(\goth g)$.
We first consider the Abelian case
and  provide a short proof of  \cite{JoSa, Theorem~1.3}  stating that $\Cal R(T)$ is a $G_\delta$.

\proclaim{Theorem 2.5}
Let $G=\Bbb R^d$ and $\Gamma=\Bbb Z^d$, $d\ge 1$.
The subsets $\Cal R(T)$ and $\Cal Ri(T)$ are both $G_\delta$ in $P(\Bbb R^d)$.
\endproclaim
\demo{Proof}
Let $(\widetilde X,\widetilde\mu)$ be the space of $\widetilde T$.
Denote by Aut$(\widetilde X,\widetilde\mu)$ the group of all $\widetilde\mu$-preserving invertible transformations of $\widetilde X$.
We endow it with the standard weak  topology.
Then Aut$(\widetilde X,\widetilde\mu)$ is a Polish group (see \cite{DaSi} and references therein).
Fix a norm on $\Bbb R^d$.
Denote by $\Cal S$ the unit ball in $\Bbb R^d$.
We define a map $\goth m:\Cal S\to\text{Aut}(\widetilde X,\widetilde\mu)$
by setting
$\goth m(v):=
\widetilde T_{v}.
$
It is obviously continuous.
We recall that the subset $\goth R$ of conservative infinite measure preserving  transformations of $(\widetilde X,\widetilde\mu)$ is a $G_\delta$  in $\text{Aut}(\widetilde X,\widetilde\mu)$ \cite{DaSi}.
It follows from this fact and Lemma~2.4(i) that the set
$$
\goth m^{-1}(\goth R)=\{v\in\Cal S\mid \text{the flow }(\widetilde T_{tv})_{t\in\Bbb R}\text{ is conservative}\}
$$
is a $G_\delta$ in  $\Cal S$, i.e. the intersection of countably many open subsets.
Since $\goth m^{-1}(\goth R)$ is centrally symmetric (i.e. if $v\in \goth m^{-1}(\goth R)$ then
$-v\in \goth m^{-1}(\goth R)$), we may assume without loss of generality that these open sets are also centrally symmetric.
The natural projection of $\Cal S$ onto $P(\Bbb R^d)$ is just the `gluing' the pairs of centrally symmetric points.
We note that the projection of $\goth m^{-1}(\goth R)$  to $P(\Bbb R^d)$ is exactly $\Cal R(\widetilde T)$.
It follows that $\Cal R(\widetilde T)$ is a $G_\delta$ in $P(\goth g)$.
It remains to apply Theorem~2.1.

To show that  $\Cal Ri(T)$ is a $G_\delta$ argue in a similar way and 
use the fact that the set of all rigid transformations  is a $G_\delta$ in Aut$(\widetilde X,\widetilde\mu)$ \cite{DaSi} and apply Lemma~2.4(iii).
\qed
\enddemo

We now consider the general case (independently of Theorem~2.5).

\proclaim{Theorem 2.6}
 The subsets $\Cal R(T)$  and $\Cal Ri(T)$  are both  $G_\delta$ in $P(\goth g)$.
\endproclaim
\demo{Proof}
Let  $\Gamma\setminus\{1\}=\{\gamma_k\mid k\in\Bbb N\}$.

(A) We first prove that  $\Cal R(T)$ is a $G_\delta$.
For each $g\in G$, the map
$$
P(\goth g)\ni\theta\mapsto \text{dist}(g,\exp(\theta)):=\inf_{h\in\exp(\theta)} \text{dist}(g,h)\in\Bbb R
\tag2-4
$$
is continuous.
Now for
 a subset $A\subset X$ with $0<\mu(A)<\infty$ and $\epsilon>0$,
we  construct  a sequence  $A_1,A_2,\dots$  of subsets in $A$ as follows (cf. with the proof of Lemma~1.5):
$$
A_1:=
\cases
A\cap T_{\gamma_1}^{-1}A, &\text{if $ \text{dist}(\gamma_1, \exp(\theta))<\epsilon$}\\
\emptyset, &\text{otherwise},
\endcases
$$
$$
A_2:=
\cases
(A\setminus(A_1\cup T_{\gamma_1}A_1))\cap T_{\gamma_2}^{-1}(A\setminus(A_1\cup  T_{\gamma_1}A_1)), &\text{if $ \text{dist}(\gamma_2, \exp(\theta))<\epsilon$}\\
\emptyset, &\text{otherwise},
\endcases
$$
and so on.
 Then (as in Lemma~1.5) $A_i\cap A_j=\emptyset$, $T_{\gamma_i}A_i\subset A$ and   $T_{\gamma_i}A_i\cap T_{\gamma_j}A_j=\emptyset $ if $i\ne j$.
For each $m\in\Bbb N$, we  set
$$
\Theta_{\epsilon,A,m}:=\bigg\{\theta\in P(\goth g)\mid \sum_{j\le m}\mu(A_j)>0.4\mu(A)\bigg\}.
$$
We note that for each $j>0$, the map $P(\goth g)\ni\theta\mapsto\mu(A_j)\in\Bbb R$ is lower semicontinuous.
Indeed, this map is (up to a multiplicative constant) is the indicator function of the subset $\{\theta\mid  \text{dist}(\gamma_j,\exp(\theta))<\epsilon\}$ which is open because \thetag{2-4} is continuous.
It follows that $\Theta_{\epsilon,A,m}$ is an open subset in $P(\goth g)$.
Fix  a countable family  $\goth D$ of subsets of finite positive measure in $X$ such that $\goth D$ is dense in $\goth B$.
We claim that
$$
\Cal R(T)=\bigcap_{D\in\goth D}\bigcap_{l=1}^\infty\bigcup_{m=1}^\infty\Theta_{1/l,D,m}.
\tag2-5
$$
Indeed, if $T$ is recurrent along a line $\theta\in P(\goth g)$ then  for each $\epsilon>0$ and each subset $A$ of positive measure,
$\mu(\bigsqcup_{j}A_j)\ge 0.5 \mu(A)$ (as in Lemma~1.5).
We then obtain that  there exists $m>0$ with $\mu(\bigsqcup_{j=1}^mA_i)>0.4\mu(A)$.
Hence $\theta\in\Theta_{\epsilon,A,m}$.
Let now $A$ run $\goth D$ and let $\epsilon$ run $\{1/l\mid  l\in\Bbb N\}$.
Then $\theta$ belongs to the right-hand side of~\thetag{2-5}.

Conversely, take $\theta$ from the right-hand side of \thetag{2-5}.
Let $A$ be a subset of $X$ of positive measure.
Then there is $D\in\goth D$ such that $\mu(A\cap D)>0.999\mu(D)$.
Take $l\in\Bbb N$.
Select $m>0$
such that $\theta\in  \Theta_{1/l,D,m}$.
Then
$$
\mu\bigg(\bigsqcup _{j\le m}D_{j}\bigg)>0.4\mu(D)
 \quad\text{ and hence }\quad
 \mu\bigg(\bigsqcup _{j\le m}T_{\gamma_j}D_{j}\bigg)> 0.4\mu(D).
 $$
Therefore there is $j<d$ with
 $\mu(T_{\gamma_j} A\cap A)>0$ and (because $\theta\in \Theta_{1/l,D,m}$) $\text{dist}(\gamma_j,\exp(\theta))<1/m$.

(B) To show that $\Cal Ri(T)$ is $G_\delta$ we first denote by $\tau$ a metric on Aut$(X,\mu)$ compatible with the weak topology.
Now it suffices to note that 
$$
\Cal Ri(T)=\bigcap_{k=1}^\infty\bigcap_{N=1}^\infty\bigcup_{\{n>N\mid\tau(T_{\gamma_n},\text{Id})<1/k\}}\{\theta\in P(\goth g)\mid \text{dist}(\gamma_n,\exp(\theta))<1/k\}.
$$ 
and use \thetag{2-4}.
\qed
\enddemo

\head
3.  $(C,F)$-construction~and~directional~recurrence of rank-one actions
\endhead

We first remind a $(C,F)$-construction of  group actions (see \cite{Da1} for a detailed exposition and various applications).

Let $(C_n)_{n>0}$ and $(F_n)_{n\ge 0}$ be two sequences of finite subsets in $\Gamma$ such that the following conditions hold:
\roster
\item"(I)"
$F_0=\{1\}$, $1\in C_n$ and $\# C_n>1$ for all $n$,
\item"(II)"
$F_nC_{n+1}\subset F_{n+1}$ for all $n$,
\item"(III)"
$F_nc\cap F_nc'=\emptyset$ for all $c\ne c'\in C_{n+1}$ and $n$ and
\item"(IV)"
$\gamma F_nC_{n+1}C_{n+2}\cdots C_m\subset F_{m+1}$ eventually in $m$  for each $\gamma\in\Gamma$ and every $n$.
\endroster
Then the infinite product space $X_n:=F_n\times C_{n+1}\times C_{n+1}\times\cdots$ is a (compact) Cantor set.
It follows from (II) and (III) that the map
$$
X_n\ni(f_n,c_{n+1},c_{n+2},c_{n+3},\dots)\mapsto(f_nc_{n+1},c_{n+2},c_{n+3},\dots)\in X_{n+1}
$$
is  a continuous embedding.
Denote by $X$ the (topological) inductive limit of the sequence $X_1\subset X_2\subset\cdots$.
Then $X$ is a locally compact Cantor set.
For a subset $A\subset F_n$, we let $[A]_n:=\{x=(f_n,c_{n+1},\dots)\in X_n\mid f_n\in A\}$.
Then $[A]_n$ is a compact open subset of $X$.
We call it an {\it $n$-cylinder}.
The family of all cylinders, i.e. the family of all compact open subsets  of $X$ is a base of the topology in $X$.
Given $\gamma\in\Gamma$ and $x\in X$, in view of (II) and (IV), there is $n$ such that
$x=(f_n,c_{n+1},\cdots)\in X_n$ and $\gamma f_n\in F_n$.
Then we let $T_\gamma x:=(\gamma f_n,c_{n+1},\dots)\in X_{n}\subset X$.
It is standard to verify that $T_\gamma$ is a well defined homeomorphism of $X$.
Moreover, $T_\gamma T_{\gamma'}=T_{\gamma\gamma'}$ for all $\gamma,\gamma'\in\Gamma$, i.e. $T:=(T_\gamma)_{\gamma\in\Gamma}$ is a continuous  action of $\Gamma$ on $X$.
It is called the {\it $(C,F)$-action of $\Gamma$ associated with  $(C_n,F_{n-1})_{n>0}$} (see \cite{dJ}, \cite{Da1}, \cite{Da3}).
This action is free and minimal.
There is a unique (up to scaling) $T$-invariant $\sigma$-finite Borel measure $\mu$ on $X$.
It is easy to compute that
$$
\mu([A]_n)=\frac{\# A}{\# C_1\cdots\# C_n}
$$
 for all subsets $A\subset F_n$, $n> 0$, provided that $\mu(X_0)=1$.
 We note that $\mu(X)=\infty$ if and only if
 $$
\lim_{n\to\infty} \frac{\# F_n}{\# C_1\cdots\# C_n}=\infty.\tag3-1
 $$
 Of course, $(X,\mu,T)$ is an ergodic conservative dynamical system.
 It is of funny rank one (see \cite{Da1} and \cite{Da3} for the definition).
 Conversely, every funny rank-one free system appears this way, i.e. it is isomorphic to a $(C,F)$-system for an  appropriately chosen sequence $(C_n,F_{n-1})_{n\ge 1}$.
 We state  without proof a lemma from \cite{Da3}.

 \proclaim{Lemma 3.1}
 Let $A$ be a finite subset  $F_n$ and let $g\in G$.
 Then $[A]_n\cap T_g[A]_n\ne\emptyset$ if and only if
 $g\in \bigcup_{m>n}AC_{n+1}\cdots C_{m}C_m^{-1}\cdots C_{n+1}^{-1}A^{-1}$.
 Furthermore, if we let
 $$
 \Cal N^{g,A}_m:=\{(a,c_{n+1},\dots,c_m)\in A\times C_{n+1}\times\cdots\times C_m\mid gac_{n+1}\cdots c_m\in AC_{n+1}\cdots C_m\}
 $$
then
 $
 \mu([A]_n\cap T_g[A]_n)
 =\lim_{m\to\infty}\frac
 {\#\Cal N^{g,A}_m}
 {\#C_1\cdots\# C_m}.
 $
 \endproclaim

 To state the next assertion we need  more notation.
 Denote by $\pi:\goth g\setminus\{0\}\to P(\goth g)$ the natural projection.
Let $\kappa$ stand for a metric on $P(\goth g)$ compatible with the topology.
Given two sequences  $(A_n)_{n=1}^\infty$ and $(B_n)_{n=1}^\infty$ of
finite subsets in $G$, we write $A_n\gg B_n$ as $n\to\infty$ if
$$
\lim_{n\to\infty}\max_{a\in A_n,b\in B_n}\kappa(\pi(\log(ab),\pi(\log(a))=0.
$$

\proclaim{Proposition 3.2}
Let $T=(T_\gamma)_{\gamma\in\Gamma}$ be a $(C,F)$-action of $\Gamma$ associated with a sequence $(C_n, F_{n-1})_{n=1}^\infty$ satisfying \rom{(I)--(IV)}.
Then
\roster
\item"\rom{(i)}"
 $\Cal R(T)\subset\bigcap_{\gamma\in\Gamma}\gamma\cdot\bigg(\bigcap_{n=1}^\infty \overline{\bigcup_{m\ge n}\pi(\log(C_n\cdots C_mC_m^{-1}\cdots C_n^{-1}\setminus\{1\}))}\bigg)$.
\item"\rom{(ii)}"
If, moreover, the group generated by all $C_j$, $j>0$, is commutative and  $C_j\setminus\{1\}\gg C_1\cdots C_{j-1}$  as $j\to\infty$ then
$$
\Cal R(T)\subset\bigcap_{\gamma\in\Gamma}\gamma\cdot\bigg(\bigcap_{n=1}^\infty \overline{\bigcup_{m\ge n}\pi(\log(C_mC_m^{-1}\setminus\{1\}))}\bigg).
$$
\item"\rom{(iii)}"
If, in addition, there is $c_j\in \Gamma$ such that $C_j=\{1,c_j\}$  for  each $j>0$ then
$$
\Cal R(T)\subset \bigcap_{\gamma\in\Gamma}\gamma\cdot\bigg(\bigcap_{n=1}^\infty
 \overline{\{\pi(\log c_m)\mid m\ge n\}}\bigg).
 $$
\endroster
\endproclaim
\demo{Proof}
(i)
Let $\theta\in\Cal R(T)$.
Then for each $n>0$,
there is a sequence $(\gamma_m)_{m=1}^\infty$ of elements of $\Gamma$ such that $\gamma_m\ne 1$ and 
$\mu(T_{\gamma_m}[1]_n\cap [1]_n)>0$ for each $m$ and $\text{dist}(\gamma_m,\exp(\theta))\to 0$ as $m\to\infty$.
Hence we deduce from  Lemma~3.1 that
$$
\inf\bigg\{\text{dist}(\gamma,\exp(\theta))\mid \gamma \in \bigcup_{m> n}C_{n+1}\cdots C_mC_m^{-1}\cdots C_{n+1}^{-1}\setminus\{1\}\bigg\}=0.
$$
This yields that
$
\theta  \in\overline{\pi\bigg(\log\bigg(\bigcup_{m> n}C_{n+1}\cdots C_mC_m^{-1}\cdots C_{n+1}^{-1}\setminus\{1\}\bigg)\bigg)}
$.
Therefore
$$
\Cal R(T)\subset\bigcap_{n\ge 1} \overline{\bigcup_{m> n}\pi(\log(C_{n+1}\cdots C_mC_m^{-1}\cdots C_{n+1}^{-1}\setminus\{1\}))}
.
$$
Since  $\Cal R(T)$ is invariant under $\Gamma$ in view of Remark~1.4(ii), the claim (i) follows.

(ii) Denote by $A$ the smallest closed  Lie subgroup of $G$ containing  all $C_j$, $j>0$.
Since $A$ is Abelian, the restriction of $\log$ to $A$ is a group homomorphism.
Hence the condition $C_j\setminus\{1\}\gg C_1\cdots C_{j-1}$ as $j\to\infty$ implies $C_jC_j^{-1}\setminus\{1\}\gg C_1C_1^{-1}\cdots C_{j-1}C_{j-1}^{-1}$ as $j\to\infty$.
Now (ii) easily follows from (i).

(iii)
It suffices to note that $C_mC_m^{-1}\setminus\{1\}=\{c_m,c_m^{-1}\}$ and $\pi(\log c_m)=\pi(\log c_m^{-1})$.

\qed
\enddemo

\head{4. Directional recurrence sets for actions of Abelian lattices}
\endhead

In this section we consider the case of Abelian $G$ in more detail.
Our purpose here is to realize various $G_\delta$-subsets of $P(\goth g)$
as $\Cal R(T)$ for rank-one actions  $T$ of $G$.
Since $G$ is simply connected, there is $d>0$ such that
$G=\Bbb R^d$.
Hence $\goth g=\Bbb R^d$ and the maps $\exp$ and $\log$ are the identities.
Replacing $\Gamma$ with an automorphic lattice we may assume without loss of generality that
$\Gamma=\Bbb Z^d$.
In the sequel we assume that $d>1$ (the case $d=1$ is trivial).
By dist$(.,.)$ we  denote the usual distance between a point and a closed subset of $\Bbb R^d$.
We also note that  $\Cal E\Cal R(T)
=\Cal R(T)$ for each measure preserving action $T$ of $\Gamma$.
We now restate Proposition~3.2 for the Abelian case.

\proclaim{Proposition 4.1}
Let $T=(T_\gamma)_{\gamma\in\Bbb Z^d}$ be a $(C,F)$-action of $\Bbb Z^d$ associated with a sequence $(C_n, F_{n-1})_{n=1}^\infty$ satisfying \rom{(I)--(IV)}.
Then
\roster
\item"\rom{(i)}"
 $\Cal R(T)\subset\bigcap_{n=1}^\infty \overline{\pi(\sum_{j\ge n}(C_j-C_j)\setminus\{0\})}$.
\item"\rom{(ii)}"
If, moreover,  $C_j\setminus\{0\}\gg C_1\cup\cdots\cup C_{j-1}$ as  $j\to\infty$ then
$$
\Cal R(T)\subset\bigcap_{n=1}^\infty \overline{\bigcup_{m\ge n}\pi((C_m-C_m)\setminus\{0\})}.
$$
\item"(iii)"
In, in addition, there is $c_j\in\Bbb Z^d$ such that $C_j=\{0,c_j\}$  for each $j>0$ then
$$
\Cal R(T)\subset \bigcap_{n=1}^\infty
 \overline{\{\pi(c_m)\mid m\ge n\}}.
 $$
\endroster
\endproclaim

The following two theorems are the main results of this section.

\proclaim{Theorem 4.2}
Let $\Delta$ be a $G_\delta$-subset of $P(\Bbb R^d)$ and let $D$  be a  countable subset of $ \Delta$.
Then there is a rank-one free infinite measure preserving action $T$ of $\Bbb Z^d$ such that $D\subset\Cal R(T)\subset \Delta$.
In particular, each countable  $G_\delta$-subset (e.g.  each  countable compact) of $P(\Bbb R^d)$ is realizable as $\Cal R(T)$ for some  rank-one free action $T$ of $\Bbb Z^d$.
\endproclaim

\demo{Proof}
Suppose first that $\Delta\ne\emptyset$.
Then without loss of generality we may think that  $D\ne\emptyset$.
Let $(\delta_n)_{n=1}^\infty$ be a sequence such that  $\delta_n\in D$ for each $n$ and  every element of $D$ occurs in this sequence infinitely many times.
Let $(\epsilon_n)_{n=1}^\infty$ be  a decreasing sequence of positive reals with $\lim_{n\to\infty}\epsilon_n= 0$.
There exists an increasing sequence $L_1\subset L_2\subset\cdots$  of closed subsets in $P(\Bbb R^d)$ such that $P(\Bbb R^d)\setminus\Delta=\bigcup_{j\ge 1}L_j$.
Let $L_1^+\subset L_2^+\subset\cdots$ be a sequence of open subsets in $P(\Bbb R^d)$ such that $ L_j^+\supset L_j$ and $\delta_j\not\in\overline{L_j^+}$ for each $j$ and $\bigcup_{j\ge 1}L_j^+\ne P(\Bbb R^d)$.
We will construct inductively two sequences $(F_n)_{n=0}^\infty$ and $(C_n)_{n=1}^\infty$ satisfying (I)--(IV) and~\thetag{3-1}.
We note in advance   that in our construction $\# C_n=2$ and $F_n$ is a symmetric cube in $\Bbb Z^d$, i.e.
there is $a_n\in\Bbb N$ such that
$$
F_n=\{(i_1,\dots,i_d)\mid -a_n< i_j\le a_n, j=1,\dots,d\},
$$
 for each $n$.
Suppose that we have defined the subsets $C_1,F_1,\dots,C_{n-1}, F_{n-1}$.
Our purpose is to construct $C_n$ and $F_n$.
Choose $c_n\in\Bbb Z^d$ such that $(c_n+F_{n-1})\cap F_{n-1}=\emptyset$,
$
\text{dist}(c_n,\delta_n)<\epsilon_n$ and
$$
\align
&\max_{f\in F_{n-1}}\text{dist}(c_n,c_n+f)<\epsilon_n,  \tag4-1\\
  &\pi(c_n)\not\in L_n^+.\tag4-2
  \endalign
$$
For that  use the fact that $\delta_n\not\in \overline{L_n^+}$.
We now let $C_n:=\{0, c_n\}$ and define $F_n$ to be a huge symmetric cube in $\Bbb Z^d$ that contains $F_{n-1}+C_n$.
Continuing this construction procedure  infinitely many times we obtain infinite sequences $(F_n)_{n=0}^\infty$ and $(C_n)_{n=1}^\infty$.
It is easy to see that (I)--(IV) and \thetag{3-1} are all satisfied.
Let $T=(T_\gamma)_{\gamma\in\Bbb Z^d}$ denote the associated $(C,F)$-action.
It is free and of rank one.
Let $(X,\mu)$ be the space of this action.

We first show that
$D\subset\Cal R(T)$.
Take $\delta\in D$, $\epsilon>0$ and a cylinder $B\subset X$.
Then there are infinitely many  $n>0$ such that  $\delta=\delta_n$
and hence  $\text{dist}(c_n,\delta)<\epsilon_n<\epsilon$.
If $n$ is large enough, $B=[B_{n-1}]_{n-1}$ for some subset $B_{n-1}\subset F_{n-1}$.
Since $[B_{n-1}]_n\subset[B_{n-1}]_{n-1}$ and $T_{c_n}[B_{n-1}]_n=[c_n+B_{n-1}]_n\subset[B_{n-1}]_{n-1}$
with $\mu([B_{n-1}]_n)=0.5\mu([B_{n-1}]_{n-1})$,
 we have
$$
\mu(T_{c_n}B\cap B)\ge\mu(T_{c_n}[B_{n-1}]_n\cap [B]_{n-1})=\mu([B_{n-1}+c_n]_n)=0.5\mu(B).
$$
Since each subset of finite measure in $X$ can be approximated with a cylinder up to an arbitrary positive real, we deduce that $\delta\in\Cal R(T)$.

We now show that  $\Cal R(T)\subset\Delta$.
It follows from  \thetag{4-1} that $\{c_n\}\gg F_{n-1}$ as $n\to\infty$.
Hence by Proposition~4.1(iii), $\Cal R(T)\subset\bigcap_{n=1}^\infty\overline{\{\pi(c_m)\mid m\ge n\}}$.
Applying \thetag{4-2}, we obtain that
$
\pi(c_m)\not\in L_m^+\supset L_n^+\supset L_n
$
for each $m\ge n$.
Hence $\Cal R(T)\cap L_n=\emptyset$ for each $n$, which yields $\Cal R(T)\subset\Delta$.

It remains to  consider the case where $\Delta=\emptyset$.
Fix $\theta\in P(\Bbb R^d)$.
Suppose that we have defined the subsets $C_1,F_1,\dots,C_{n-1}, F_{n-1}$.
Choose $c_n\in\Bbb Z^d$ such that $(c_n+F_{n-1})\cap F_{n-1}=\emptyset$,
\thetag{4-1} is satisfied,
 $$
 \align
 &\pi(c_{n})
 \text{ is up to $\epsilon_n$ close to $\theta$ (in the metric on $P(\Bbb R^d)$)
 and}
 \tag4-3 \\
& \min_{f\in F_{n-1}-F_{n-1}} \text{dist}(c_n+f,\theta)>10.\tag4-4
  \endalign
$$
We now let $C_n:=\{0, c_n\}$ and define $F_n$ to be a huge symmetric cube in $\Bbb Z^d$ that contains $F_{n-1}+C_n$.
Continuing infinitely many times we obtain infinite sequences $(F_n)_{n=0}^\infty$ and $(C_n)_{n=1}^\infty$.
It is easy to see that (I)--(IV) and \thetag{3-1} are all satisfied.
Let $T=(T_\gamma)_{\gamma\in\Bbb Z^d}$ denote the associated $(C,F)$-action.
It follows from Proposition~4.1(iii), \thetag{4-1} and \thetag{4-3} that  $\Cal R(T)\subset \{\theta\}$.
If $T$ were recurrent along $\theta$ then there is $\gamma\in\Bbb Z^d$ such that  $\gamma\ne 0$, $\text{dist}(\gamma,\theta)<0.1$ and $\mu([0]_n\cap T_\gamma[0]_n)>0$.
It follows from Lemma~3.1 that there is $l>n$ such that, $\gamma\in F_{l-1}-F_{l-1}+c_l$.
 This contradicts to~\thetag{4-4}.
  Thus we obtain that $\Cal R(T)=\emptyset$.
\qed

\enddemo

\proclaim{Theorem~4.3}
There is  a rank-one free infinite measure preserving action $T$ of $\Bbb Z^d$ such that $\Cal R(T)=P(\Bbb R^d)$.
\endproclaim
\demo{Proof}
Given $t \in\Bbb N$ and $N>0$, we let
$$
\Cal K_{t,N}:=\{(i_1,\dots,i_d)\in\Bbb Z^d\mid  |i_j|<N\text{ and $t$ divides $i_j$, $j=1,\dots,d$}\}.
$$
Then for each $\epsilon>0$ and each integer $t>0$, there is $N>0$ such that
$$
\sup_{\delta\in P(\Bbb R^m)}\min_{0\ne \gamma\in\Cal K_{t,N}}\text{dist}(\gamma,\delta)<\epsilon.
\tag4-5
$$
Fix a sequence of positive reals $\epsilon_n$, $n\in\Bbb N$, decreasing to $0$.
We will construct inductively the sequences $(F_{n-1})_{n>0}$ and $(C_{n})_{n>0}$
satisfying (I)--(IV) and \thetag{3-1}.
As usual, $F_0=\{0\}$.
Suppose we have defined $(F_{j}, C_j)_{j=1}^n$.
Suppose that $F_n$ is a symmetric cube.
Denote by $t_n$ the length of an edge of this cube.
We now construct $C_{n+1}$ and $F_{n+1}$.
By~\thetag{4-5}, there is $N_n$ such that
$\min_{0\ne \gamma\in\Cal K_{3t_n,N_n}}\text{dist}(\gamma,\delta)<\epsilon_n$
for each $\delta\in P(\Bbb R^d)$.
 Let $C_{n+1}:=\Cal K_{3t_n,M_n}$, where $M_n$ is an integer large so that
 $$
 \#\{\gamma\in \Cal K_{3t_n,M_n}\mid \gamma+\Cal K_{3t_n,N_n}\subset \Cal K_{3t_n,M_n}\}>0.5\#\Cal K_{3t_n,M_n}.\tag4-6
 $$
Now  let $F_{n+1}$ be a huge symmetric cube in $\Bbb Z^d$ such that $F_{n+1}\supset F_n+C_{n+1}$.
  Continuing  this construction process infinitely many times we define the infinite sequences
 $(F_n)_{n\ge 0}$ and $(C_n)_{n\ge 1}$ as desired.
 Let $T$ be the $(C,F)$-action of $\Bbb Z^d$ associated with these sequences.
 It is  free and of rank-one.
 Denote by  $(X,\mu)$ the space of this action.
  We claim that  $\Cal R(T)=P(\Bbb R^d)$.
  Indeed, take $\epsilon>0$, $\delta\in P(\Bbb R^d)$ and a cylinder $B\subset X$.
 Then there is $n>0$ and a subset $B_n\subset F_n$ such that $B=[B_n]_n$ and $\epsilon_n<\epsilon$.
 There is $\gamma\in\Cal K_{3t_n,N_n}\setminus \{0\}$ such that $\text{dist}(\gamma,
\delta)<\epsilon_n$.
By~\thetag{4-6},
 $\#(C_{n+1}\cap(C_{n+1}-\gamma))\ge 0.5\# C_{n+1}$.
 Therefore
$$
\align
\mu(T_\gamma B\cap B)&\ge\mu(T_\gamma[B_n+(C_{n+1}\cap(C_{n+1}-\gamma))]_{n+1}\cap [B_n]_n)\\
&=\mu([B_n+(C_{n+1}\cap(C_{n+1}+\gamma))]_{n+1})\\
&\ge 0.5\mu(B).
\endalign
$$
 The standard approximation argument implies that $T$ is recurrent along $\delta$.
 \qed
\enddemo

\remark{\rom{Remark 4.4}}
\roster
\item"\rom{(i)}"
If we choose $M_m$ in the above construction large so that the  inequality
$$
 \#\{\gamma\in \Cal K_{3t_n,M_n}\mid \gamma+\Cal K_{3t_n,N_n}\subset \Cal K_{3t_n,M_n}\}>(1-n^{-1})\#\Cal K_{3t_n,M_n}.
 $$
 holds in place of \thetag{4-6} then the corresponding $(C,F)$-action $T$ will possess the stronger property
 $\Cal Ri(T)=P(\Bbb R^d)$.
 \item"\rom{(ii)}"
 In a similar way, the statement of Theorem~4.2 remains true if we replace $\Cal R(T)$ with $\Cal Ri(T)$.
 \endroster
 \endremark

\head 5. Generic $\Bbb Z^d$-action is recurrent in every direction
 \endhead

 Let $(X,\mu)$ be a $\sigma$-finite non-atomic standard measure space.
 We recall that  the group of all $\mu$-preserving invertible transformations of $X$
 is denoted by  Aut$(X,\mu)$.
 It is  endowed with the weak (Polish) topology.
 For a nilpotent Lie group  $G$, we
 denote by $\Cal A^G_\mu$ the set of all $\mu$-preserving actions of $G$ on $(X,\mu)$.
 We consider every element $A\in \Cal A^G_\mu$ as a continuous homomorphism $g\mapsto A_g$ from $G$ to Aut$(X,\mu)$.
 The group  Aut$(X,\mu)$ acts  on $\Cal A^G_\mu$ by conjugation, i.e. $(S\cdot A)_g:=SA_gS^{-1}$ for all $g\in G$, $S\in\text{Aut}(X,\mu)$ and $A\in \Cal A^G_\mu$.
 We endow $\Cal A^G_\mu$ with the compact-open topology, i.e. the topology of uniform convergence on the compact subsets of $G$.
 
The following lemma is well known.
We state it without proof.

 \proclaim{Lemma 5.1} $ \Cal A^G_\mu$
 is a Polish space.
 The action of \rom{Aut}$(X,\mu)$ on this space is continuous.
 \endproclaim

 Let $S^1$ be the unit sphere in $\goth g$ and let $K:=\exp (S^1)$.

 \proclaim{Lemma 5.2} Let $\mu(X)=1$.
 Then the subset
 $$
\Cal Z:= \{A\in  \Cal A^G_\mu\mid h(A_{g}) =0\text{ for each }g\in K\}
 $$
 is an invariant
 $G_\delta$ in $ \Cal A^G_\mu$.
 \endproclaim

 \demo{Proof} Denote by $\Cal P$ the set of all finite partitions of $X$.
 Fix a countable subset $\Cal P_0\subset\Cal P$ which is dense in $\Cal P$ in the natural topology.
 For each $P\in\Cal P_0$ and $n>0$, the map
 $$
  \Cal A^G_\mu\times K\ni (A,g)\mapsto H\bigg(P\,\bigg|\,\bigvee_{j=1}^n A^{-j}_gP\bigg)\in\Bbb R
 $$
 is continuous.
 Therefore the map
 $$
m_{P,n}: \Cal A^G_\mu\ni A\mapsto m_{P,n}(A):=\max_{g\in K}H\bigg(P\,\bigg|\,\bigvee_{j=1}^n A^{-j}_gP\bigg)\in\Bbb R
 $$
 is well defined and continuous.
 Hence the subset
 $$
 \Cal Z':=\bigcap_{P\in\Cal P_0}\bigcap_{r=1}^\infty\bigcap_{N=1}^\infty\bigcup_{l>N}\bigg\{A\in \Cal A^G_\mu\,\bigg|\, m_{P,l}(A)<1/r\bigg\}
 $$
 is a $G_\delta$ in $ \Cal A^G_\mu$.
 We now show that $\Cal Z'=\Cal Z$.
 It is easy to see that  $\Cal Z'\subset\Cal Z$ because $h(A_g)=\sup_{P\in\Cal P_0}H(P\mid\bigvee_{j=1}^\infty A_g^{-j}P)$.
 Conversely, let $A\in\Cal Z$.
Fix $P\in\Cal P_0$, $r>1$ and $N>0$.
Then for each $g\in K$, there is $l_g>N$ such that  $H(P\mid\bigvee_{j=1}^{l_g}A_g^{-j}P)<1/r$.
Of course, this inequality holds in a neighborhood of $g$ in $G$.
Since $K$ is compact and the map $\Bbb N\ni n\mapsto H(P\mid\bigvee_{j=1}^{n}A_g^{-j}P)$ decreases, there is $l>N$ such that
$H(P\mid\bigvee_{j=1}^lA_g^{-j}P)<1/r$ for all $g\in K$, i.e. $m_{P,l}(A)<1/r$.
 This means that $A\in\Cal Z'$.
It is obvious  that $\Cal Z$ is Aut$(X,\mu)$-invariant.
 \qed
 \enddemo

 Let $\Gamma$ be a co-compact lattice in $G$.
 Fix a a cross-section $s:G/\Gamma\to G$ of the natural projection $G\to G/\Gamma
 $ such that the subset $s(G/\Gamma)$ is relatively compact in $G$.
 Denote by $h_s$ the corresponding 1-cocycle.
 Given a $\Gamma$-action $T$ on $(X,\mu)$, we construct (via $h_s$) the induced $G$-action   $\widetilde T$ on the space
 $(G/\Gamma\times X,\lambda\times\mu)$.
  In the following lemma we show that the ``inducing'' functor is continuous.

 \proclaim{Lemma 5.3}
The map $\Cal A^\Gamma_\mu\ni T\mapsto\widetilde T\in\Cal A^G_{\lambda_{G/\Gamma}\times\mu}$ is continuous.
 \endproclaim

 \demo{Idea of the proof}
  It is enough to note that for each compact subset $K\subset G$, the set
 $F:=\{h_s(g,y)\mid g\in K, y\in G/\Gamma\}\subset\Gamma$ is finite.
 Therefore, given two $\Gamma$-actions $T$ and $T'$, if the transformation $T_\gamma$ is ``close'' to $T'_\gamma$ for each $\gamma\in F$ then the transformation
 $\widetilde T_g$ is ``close'' to $\widetilde T'_g$ uniformly on $K$.
 \qed
 \enddemo

 From now on let $\mu(X)=\infty$.
 Denote by $(X^\bullet,\mu^\bullet)$ the Poisson suspension of $(X,\mu)$.
 Given $R\in\text{Aut}(X,\mu)$, let $R^\bullet$ stand for  the Poisson suspension of $R$ (see \cite{Ro}, \cite{Ja--Ru}).
 We note that $\text{\rom{Aut}}(X^\bullet,\mu^\bullet)$ is a topological
 $\text{\rom{Aut}}(X,\mu)$-module.

 \proclaim{Lemma 5.4}
 The map \rom{Aut}$(X,\mu)\ni R\mapsto R^\bullet\in\text{\rom{Aut}}(X^\bullet,\mu^\bullet)$ is a continuous homomorphism.
 \endproclaim
\demo{Idea of the proof}
Let $U_R$ and $U_{R^\bullet}$ denote the Koopman unitary operators generated by $R$ and $R^\bullet$ respectively.
Then it is enough to note that $U_{R^\bullet}$ is unitarily equivalent in a canonical way to the exponent $\bigoplus_{n\ge 0}U_R^{\odot n}$ (see \cite{Ne}, \cite{Ro}) and the map $U_R\mapsto U_R^{\odot n}$ is continuous in the weak operator topology for each $n$.
\qed
\enddemo

\proclaim{Lemma 5.5}
Let a transformation $R\in\text{\rom{Aut}}(X,\mu)$ be non-conservative.
If there is an ergodic  countable transformation subgroup $N\subset \text{\rom{Aut}}(X,\mu)$ such that
$$
\{SR^nx\mid n\in\Bbb Z\}=\{R^nSx\mid n\in \Bbb Z\}\quad\text{at a.e. }x\in X\text{ for each }S\in N\tag5-1
$$
then $R^\bullet$ is a Bernoulli transformation of infinite entropy.
\endproclaim

\demo{Proof}
We consider Hopf decomposition of $X$, i.e. a partition of $X$ into two $ R$-invariant subsets $X_d$ and $X_c$ such that the restriction of $R$ to $X_d$ is totally dissipative and the restriction of $R$ to $X_d$ is conservative (see \cite{Aa}).
By condition of the lemma, $\mu(X_d)>0$.
It follows from \thetag{5-1} that $X_d$ is invariant under $N$.
Since $N$ is ergodic, $\mu(X_c)=0$, i.e. $R$ is totally dissipative, i.e. there is a subset $W\subset X$ such that $X=\bigcup_{n\in\Bbb Z}R^nW$ (mod 0) and $R^nW\cap T^mW=\emptyset$ if $n\ne m$.
Therefore $R^\bullet$ is Bernoulli \cite{Ro}.
Since $\mu\restriction W$ is not purely atomic,
 $h(R^\bullet)=\infty$  \cite{Ro}.
\qed
\enddemo

 We now state  the main result of this section.

 \proclaim{Theorem 5.6}
 The subset $\Cal V$ of $\Bbb Z^d$-actions $T$ on $(X,\mu)$ with $\Cal R(T)=P(\Bbb R^d)$ is residual in $\Cal A_\mu^{\Bbb Z^d}$.
 \endproclaim

\demo{Proof}
Let $\lambda$ denote the Lebesgue measure on the torus $\Bbb R^d/\Bbb Z^d$.
If follows from Lemmata~5.3 and 5.4 that the mapping
$$
\Cal A^{\Bbb Z^d}_\mu\ni T\mapsto \widetilde T^\bullet\in\Cal A^{\Bbb R^d}_{\lambda\times\mu}
$$
is continuous.
Let
$
\Cal Z:= \{A\in  \Cal A^{\Bbb R^d}_{(\lambda\times\mu)^\bullet}\mid h(A_{g}) =0\text{ for each }g\in \Bbb R^d\}.
 $
By Lemma~5.2, $\Cal Z$ is a $G_\delta$ in $\Cal A^{\Bbb R^d}_\mu$.
Hence the subset $\Cal W:=\{T\in \Cal A^{\Bbb Z^d}_\mu\mid \widetilde T^\bullet\in\Cal Z\}$ is an
$G_\delta$ in $\Cal A^{\Bbb Z^d}_\mu$.
Of course, $\Cal W$ is Aut$(X,\mu)$-invariant.
It is well known that the subset $
\Cal E:=\{T\in \Cal A^{\Bbb Z^d}_\mu\mid \text{$T$ is ergodic}\}$ is an
Aut$(X,\mu)$-invariant
$G_\delta$ in $\Cal A^{\Bbb Z^d}_\mu$.
Hence the intersection
$\Cal W\cap \Cal E$ is also an Aut$(X,\mu)$-invariant
$G_\delta$ in $\Cal A^{\Bbb Z^d}_\mu$.
Take an action $T\in\Cal A^{\Bbb Z^d}\cap\Cal E$ and a line $\theta\in P(\Bbb R^d)$.
If $\theta\not\in\Cal R(T)$ then $\theta\not\in\Cal R(\widetilde T)$.
Since $T$ is ergodic, $\widetilde T$ is also ergodic.
Hence the $\Bbb Q^d$-action $(\widetilde T_q)_{q\in\Bbb Q^d}$ is also ergodic.
Then by Lemma~5.5, $h(\widetilde T_r^\bullet)=\infty$ for each $r\in\theta$, $r\ne 0$.
Therefore $T\not\in\Cal W$.
This yields that $\Cal W\cap\Cal E\subset\Cal V$.
It remains to show that $\Cal W\cap\Cal E$ is dense in $\Cal A_\mu^{\Bbb Z^d}$.
Let $T$ be an ergodic  free action of $\Bbb Z^d$
such that $\Cal Ri(T)=P(\Bbb R^d)$ (see Remark~4.4(i) and Theorem~4.3).
By Theorem~2.1, $\Cal Ri(\widetilde T)=P(\Bbb R^d)$.
 Then in view of  Lemma~5.4, for each $g\in\Bbb R^d$, the transformation $\widetilde T_g^\bullet$ is rigid.
 Hence $h(\widetilde T_g^\bullet)=0$.
 Thus, $T\in\Cal W\cap\Cal E$.
 It follows from Rokhlin lemma for the infinite measure preserving free $\Bbb Z^d$-actions that the conjugacy class of $T$, i.e. the Aut$(X,\mu)$-orbit of $T$, is dense in
$\Cal A_\mu^{\Bbb Z^d}$ (see, e.g. \cite{DaSi}).
Of course, the conjugacy class of $T$ is a subset of $\Cal W\cap\Cal E$.
\qed
\enddemo

Using some ideas from the proof of the above theorem  we can show  the following proposition.

 \proclaim{Proposition 5.7} There is a Poisson action\footnote{We recall that a probability preserving action of a group $G$ is called Poisson if it is isomorphic to the Poisson suspension of an infinite measure preserving action of $G$.} $V$ of $\Bbb R^d$ of $0$ entropy such that
for each $0\ne g\in \Bbb R^m$, the transformation  $V_g$ is Bernoullian and of infinite entropy.
 \endproclaim
\demo{Proof}
By Theorem~4.2, there exists rank-one (by cubes) infinite measure preserving action $T$ of $\Bbb Z^d$
such that $\Cal R(T)=\emptyset$.
Then $\widetilde T^\bullet$ is a Poisson (finite measure preserving) action of $\Bbb R^d$.
We note that $h(\widetilde T^\bullet)=h(\widetilde T^\bullet\restriction \Bbb Z^d)=h((\widetilde T\restriction \Bbb Z^d)^\bullet)$.
We note $\widetilde T\restriction \Bbb Z^d=I\times T$, where $I$ denotes the trivial action of $\Bbb Z^d$ on the torus $(\Bbb R^d/\Bbb Z^d,\lambda)$.
It follows from \cite{Ja--Ru} that $h((I\times T)^\bullet)=h(T^\bullet)$.
Since $T$ is of rank one, $h(T^\bullet)=0$ by \cite{Ja--Ru}\footnote{This fact was proved in \cite{Ja--Ru} only for $d=1$. However in the general case the proof is similar.}.
Thus we obtain that $h(\widetilde T^\bullet)=0$.
On the other hand, arguing as in the proof of Theorem~5.6, we deduce from Theorem~2.1 and Lemma~5.5 that for each $g\in\Bbb R^d\setminus\{0\}$, the transformation
$\widetilde T^\bullet_g$ is  Bernoulli and of infinite entropy.
\qed
\enddemo

In a similar way, using Remark~4.4(ii) we can show the following more general statement.

\proclaim{Proposition 5.8}
Let $\Delta$ be a $G_\delta$-subset of $P(\Bbb R^d)$ and let $D$  be a  countable subset of $ \Delta$.
Then there is a Poisson action  $V$ of $\Bbb R^d$ of 0 entropy such that for each nonzero $g\not\in\bigcup_{\theta\in\Delta}\theta$, the transformation
$V_g$ is Bernoulli and of infinite entropy and for each $g\in \bigcup_{\theta\in D}\theta$, the transformation $V_g$ is rigid (and hence of 0 entropy).
\endproclaim

\head
6. Directional recurrence for actions of the Heisenberg group
\endhead

Consider now
 the 3-dimentional real Heisenberg group  $H_3(\Bbb R)$ which is perhaps the simplest example of a non-commutative simply connected nilpotent Lie group.
 We recall that
 $$
 H_3(\Bbb R)=\Bigg\{\pmatrix
 1 & t_1 & t_3\\
 0& 1 & t_2\\
 0& 0 & 1
 \endpmatrix \Bigg| \, t_1,t_2,t_3\in\Bbb R
 \Bigg\}.
 $$
 We introduce the following notation:
 $$
 a(t):=
\pmatrix
 1 & t & 0\\
 0& 1 & 0\\
 0& 0 & 1
 \endpmatrix ,\
 b(t):=
\pmatrix
 1 & 0 & 0\\
 0& 1 & t\\
 0& 0 & 1
 \endpmatrix ,\
 c(t):=
\pmatrix
 1 & 0 & t\\
 0& 1 & 0\\
 0& 0 & 1
 \endpmatrix
 .
  $$
  Then the maps $\Bbb R\ni t\mapsto a(t)\in H_3(\Bbb R)$, $\Bbb R\ni t\mapsto b(t)\in H_3(\Bbb R)$,$\Bbb R\ni t\mapsto c(t)\in H_3(\Bbb R)$ are continuous  homomorphisms,
  the subset $\{c(t)\mid t\in\Bbb R\}$ is the center of $H_3(\Bbb R)$, $a(t_1)b(t_2)=b(t_2)a(t_1)c(t_1t_2)$ for all $t_1,t_2\in\Bbb R$ and
  $$
 \pmatrix
 1 & t_1 & t_3\\
 0& 1 & t_2\\
 0& 0 & 1
 \endpmatrix
 =c(t_3)b(t_2)a(t_1) \quad \text{for all }t_1,t_2,t_3\in\Bbb R.
  $$
  We also note that  the  Lie algebra  of $H_3(\Bbb R)$ is
$$
\goth h_3(\Bbb R):=
\Bigg\{
\pmatrix
 0& t_1 & t_3\\
 0& 0 & t_2\\
 0& 0 & 0
 \endpmatrix\Bigg|\, \alpha,\beta,\gamma\in\Bbb R
 \Bigg\}.
$$
The exponential map $\exp:\goth h_3(\Bbb R)\to H_3(\Bbb R)$ is given by the formula
$$
\exp\pmatrix
 0& t_1 & t_3\\
 0& 0 & t_2\\
 0& 0 & 0
 \endpmatrix=
  \pmatrix
 1 & t_1 & t_3+\frac{t_1t_2}2\\
 0& 1 & t_2\\
 0& 0 & 1
 \endpmatrix.
 $$
 The adjoint action of $H_3(\Bbb R)$ on $\goth h_3(\Bbb R)$ is given by the formula
 $$
 \pmatrix
 1&x&z\\
 0&1&y\\
 0&0&1
 \endpmatrix
 \cdot
 \pmatrix
 0&\alpha&\gamma\\
 0&0&\beta\\
 0&0&0
 \endpmatrix
 =\pmatrix
 0&\alpha& \gamma+x\beta-y\alpha\\
 0&0&\beta\\
 0&0&0
 \endpmatrix
 .
 $$
  We also give an example of a right-invariant metric $d$ on $H_3(\Bbb R)$:
 $$
 d(c(t_3)b(t_2)a(t_1),c(t_3')b(t_2')a(t_1')):=|t_1-t_1'|+|t_2-t_2'|+ |t_3-t_3'+t_2'(t_1'-t_1)|.
 $$
Let  $\Gamma$ be a lattice in $H_3(\Bbb R)$.
It is well known (see, e.g. \cite{DaLe}) that  there is $k>0$ such that $\Gamma$ is automorphic to the following lattice:
$$
\{c(n_3/k)b(n_2)a(n_1)\mid n_1,n_2,n_3\in\Bbb Z\}.
$$
From now on we will assume that $k=1$ and hence
$$
\Gamma=H_3(\Bbb Z):=
\{c(n_3)b(n_2)a(n_1)\mid n_1,n_2,n_3\in\Bbb Z\}.
$$
Let $F_n:=\{c(j_3)b(j_2)a(j_1)\mid |j_1|<L_n, |j_2|<L_n,|j_3|<M_n\}$, where $L_n$ and $M_n$ are positive integers.
It is easy to verify that if  $L_n\to\infty$, $M_n\to\infty$ and $L_n/M_n\to 0$ as $n\to\infty$ then $(F_n)_{\ge 1}$ is a F{\o}lner sequence in $H_3(\Bbb Z)$.

In the following three theorems we construct rank-one actions of $H_3(\Bbb Z)$ with various sets of recurrence and rigidity: empty, countable and uncountable.

\proclaim{Theorem 6.1}
There is a rank-one free  infinite measure preserving  action $T$ of $H_3(\Bbb Z)$ such that $\Cal R(T)=\emptyset$.
\endproclaim
\demo{Proof}
Let $C_n:=\{1, a(t_n)\}$, where $(t_n)_{n\in\Bbb N}$ is a sequence of integers that grows fast,
and let $(F_n)_{n\ge 0}$ be a F{\o}lner sequence in $H_3(\Bbb R)$ such that (I)--(IV) and \thetag{3-1} are satisfied
and, in addition, $C_n \setminus \{1\} \gg C_1 \cdots C_{n-1}$ as $n\to\infty$.
Denote by $T$ the $(C,F)$-action  of $H_3(\Bbb Z)$ associated with $(C_n, F_{n-1})_{n\in\Bbb N}$.
Let $\theta\in P(\goth h_3(\Bbb R))$ stand
for the line in $\goth h_3(\Bbb R)$ passing through the vector $\pmatrix
 0&1&0\\
 0&0&0\\
 0&0&0
 \endpmatrix$.
 Since $\pi(\log a(t_n))=\theta$,
we deduce from  Proposition~3.2(iii),
$$
\Cal R(T)\subset \bigcap_{\gamma\in\Gamma}\gamma\cdot\bigg(\bigcap_{n=1}^\infty
 \overline{\{\pi(\log a(t_m))\mid m\ge n\}}\bigg)\subset
 \bigcap_{\gamma\in\Gamma}
 \{\gamma\cdot
\theta\}=\emptyset.
 $$
 \qed
\enddemo

Given $t\in\Bbb R$, let $\theta_t\in P(\goth h_3(\Bbb R))$ be the line in $\goth h_3(\Bbb R)$ passing through the vector $\pmatrix
 0&1&t\\
 0&0&0\\
 0&0&0
 \endpmatrix$.
Then $\exp(\theta_t)\ni c(t)a(1)$.
We also denote by $\theta_\infty$ the line in $\goth h_3(\Bbb R)$ passing through  the vector
$\pmatrix
 0&0&1\\
 0&0&0\\
 0&0&0
 \endpmatrix$.
Of course, the set $\{\theta_l\mid l\in\Bbb Z\}$ is the  $H_3(\Bbb Z)$-orbit $\{\gamma\cdot\theta_0\mid\gamma\in H_3(\Bbb Z\}$ of $\theta_0$.
The point
$\theta_\infty$ is the only limit point of this orbit  in $P(\goth h_3(\Bbb R))$.
In a similar way,
the set $\{\theta_t\mid t\in\Bbb R\}$ is the  $H_3(\Bbb R)$-orbit of $\theta_0$.
The closure of this orbit is the union of this orbit with the limit point 
$\theta_\infty$.

\proclaim{Theorem 6.2}
There is a rank-one free  infinite measure preserving  action $T$ of $H_3(\Bbb Z)$ such that
$\Cal R(T)=\{\theta_l\mid l\in\Bbb Z\}\cup\{\theta_\infty\}$.
Therefore $\Cal E\Cal R(T)=\{\theta_\infty\}$ and hence $\Cal R(T)\ne\Cal E\Cal R(T)$.
\endproclaim
\demo{Proof}
We let
$$
\align
F_n&:=\{c(j_3)b(j_2)a(j_1)\mid |j_1|<L_n, |j_2|<L_n,|j_3|<M_n\}\text{ and}\\
C_n&:=\{ c(ik_n)a(jk_n)\mid j=0,1\text{ and } |i|\le I_n\},
\endalign
$$
where  $(L_n)_{n\ge 1}$, $(M_n)_{n\ge 1}$, $(k_n)_{n\ge 1}$ and   $(I_n)_{n\ge 1}$  are sequence of integers  chosen in such a way such that
\roster
\item"$(\bullet)$"
(I)--(IV) from Section~3 and \thetag{3-1} are satisfied
\item"$(*)$"
 $C_n \setminus \{1\} \gg C_1 \cdots C_{n-1} $ as $n\to\infty$,
  \item"$(\diamond)$"
 $L_n\to\infty$, $M_n\to\infty$, ${L_n/M_n}\to 0$ and
 \item"$(\circ)$"
 $ I_n\to+\infty$,
$ L_{n-1}/I_n\to 0$.
\endroster
Denote by $T$ the $(C,F)$-action  of $H_3(\Bbb Z)$ associated with $(C_n, F_{n-1})_{n\in\Bbb N}$.
It is well defined in view of $(\bullet)$.
Moreover, $(F_n)_{n\ge 1}$ is a F{\o}lner sequence in $H_3(\Bbb Z)$ in view of~$(\diamond)$.
It is standard to verify that
$$
\overline{\bigcup_{m>n}\pi(\log(C_mC_m^{-1}\setminus\{1\}))}=\{\theta_l\mid l\in\Bbb Z\}
\cup\{\theta_\infty\}
$$
for each $n>0$.
Hence by Proposition~3.2(ii), $\Cal R(T)\subset\{\theta_n\mid n\in\Bbb Z\}\cup\{\theta_\infty\}$.
In view of Remark~1.4(ii), to prove the converse inclusion it suffices to show that $\theta_1,\theta_\infty\in\Cal R(T)$.
For $n\ge 1$, take a subset $D\subset F_{n-1}$.
It follows from the definition of $F_{n-1}$ that for each $\gamma\in D$, there is $j\in\Bbb Z$ such that $|j|<L_{n-1}$ and $a(k_n)\gamma a(-k_n)=\gamma c(jk_n)$.
Let
$$
C_n':=\{w\in C_{n}\mid c(jk_n)a(k_n)w\in C_n \text{ whenever }|j|<L_{n-1}\}.\tag6-1
$$
Then $C_n'=\{c(ik_n)\mid  |i|<I_n,|i\pm L_{n-1}|<I_n\}$.
Hence
$\#C_n'/\# C_n\to 1/2$ as $n\to\infty$  in view of~$(\circ)$ and hence 
$$
\max_{D
\subset F_{n-1}}
| \mu([D]_{n-1})/\mu([DC_n']_{n})- 1/2|  \to 0\tag6-2
$$
as $n\to\infty$.
On the other hand, in view of \thetag{6-1}, we have
$$
 T_{a(k_n)}[DC_n']_{n}=\bigsqcup_{\gamma\in D}T_{a(k_n)}[\gamma C_n']_n
 =\bigsqcup_{\gamma\in D}[a(k_n)\gamma a(-k_n)a(k_n)C_n']_n
 \subset \bigsqcup_{\gamma\in D}[\gamma C_n]_n.
 $$
 Thus $  T_{a(k_n)}[DC_n']_{n}\subset[D]_{n-1}.$
 Since $a(k_n)\in\exp(\theta_1)$ and \thetag{6-2} holds,
 it follows that $T$ is recurrent along $\theta_1$.
 To prove that $\theta_\infty\in\Cal R(T)$, we let
 $$
C_n'':=\{w\in C_{n}\mid c(k_n)w\in C_n \}.
$$
Then $\#C_n''/\# C_n\to 1$ and
and hence $\max_{D
\subset F_{n-1}}
| \mu([D]_{n-1})/\mu([DC_n']_{n})-1|  \to 0$
as $n\to\infty$.
Moreover,
$T_{c(k_n)}[DC_{n}'']_n\subset[DC_n]_n=[D]_{n-1}$.
Hence $T$ is recurrent along $\theta_\infty$.
\qed
\enddemo

\proclaim{Theorem 6.3}
There is a rank-one free  infinite measure preserving  action $T$ of $H_3(\Bbb Z)$ such that
$\Cal R(T)=\Cal Ri(T)=\{\theta_t\mid t\in\Bbb R\}\cup\{\theta_\infty\}=\Cal E\Cal R(T)=\Cal E\Cal Ri(T)$.
\endproclaim
\demo{Proof}
Let
$$
\align
F_n&:=\{c(j_3)b(j_2)a(j_1)\mid |j_1|<L_n, |j_2|<L_n,|j_3|<M_n\},\\
C_n &:=\{c(jk_n)a(ik_n)\mid |j|\le l_nJ_n, |i|\le l_nI_n\},\\
C_n^0 &:=\{c(jk_n)a(ik_n)\mid |j|\le l_n, |i|\le l_n\},\\
\endalign
$$
where  $(L_n)_{n\ge 1}$, $(M_n)_{n\ge 1}$, $(k_n)_{n\ge 1}$,  $(I_n)_{n\ge 1}$,  $(J_n)_{n\ge 1}$ and $(l_n)_{n\ge 1}$ are sequence of integers  such that $(\bullet)$, $(*)$, $(\diamond)$ hold,
\roster
  \item"$(\vartriangle)$"
  $\sup_{t\in\Bbb R\cup\{\infty\}}\min_{1\ne \gamma\in C_n^0}\text{dist}(\gamma,\theta_t)<1/n$
  and
 \item"$(\blacktriangle)$"
 $\#(\{w\in C_n\mid \bigcup_{d\in F_{n-1}}\bigcup_{c\in C_n^0}d^{-1}cdw\subset C_n\})>(1-1/n)\#C_n$
\endroster
for each $n\in\Bbb N$.
Denote by $T$ the $(C,F)$-action  of $H_3(\Bbb Z)$ associated with $(C_n, F_{n-1})_{n\in\Bbb N}$.
It is standard to verify that
$$
\overline{\bigcup_{m>n}\pi(\log(C_mC_m^{-1}\setminus\{1\}))}=\{\theta_t\mid t\in\Bbb R\}\cup\{\theta_\infty\}.
$$
Hence by Proposition~3.2(ii), $\Cal R(T)\subset\{\theta_t\mid t\in\Bbb R\}\cup\{\theta_\infty\}$.
To prove the converse inclusion, we take  $\theta_t$ for some $t\in\Bbb R\cup\{\infty\}$.
By $(\vartriangle)$,  there is $\gamma\in C_n^0\setminus \{1\}$ such that dist$(\gamma,\theta_t)<1/n$.
Let 
$$
C_n':=\bigg\{w\in C_{n}\,\bigg|\, \bigcup_{d\in F_{n-1}}d^{-1}\gamma d w C_n^0\subset C_n\bigg\}.
$$
Then $\#C_n'/\# C_n>1-1/n$ in view of $(\blacktriangle)$  and hence 
for
each
subset $D\subset F_{n-1}$,
we have
$\mu([D]_{n-1}\setminus[DC_n']_{n})<\mu([D]_n)/n$.
On the other hand, 
 $$
 T_{\gamma}[DC_n']_{n}=\bigsqcup_{d\in D}T_{\gamma}[dC_n']_n
 =\bigsqcup_{d\in D}[dd^{-1}\gamma d C_n']_n
 \subset \bigsqcup_{d\in D}[dC_n]_n=[D]_{n-1}.
 $$
 It follows that $T$ is rigid along $\theta_t$.
 Thus we showed that  
 $\{\theta_t\mid t\in\Bbb R\}\cup\{\theta_\infty\}\subset\Cal Ri(T)$.
 \qed
\enddemo

\head 7. Some open problems
and  concluding remarks  
\endhead

\roster

\item
Which $G_\delta$-subsets of $P(\goth g)$ are realizable as $\Cal R(T)$ or $\Cal Ri(T)$ for an ergodic infinite measure preserving action $T$ of $\Gamma$?
In particular, let $\theta\in P(\goth g)$.
Whether the subset $P(\goth g)\setminus\{\theta\}$ is realizable?
In the case where $G=\Bbb R^2$ and $\Gamma=\Bbb Z^2$, $P(\goth g)$ is homeomorphic to the circle.
Whether a proper arc of this circle is realizable? 
\item
Suppose that a subset of $P(\goth g)$ is realizable as $\Cal R(T)$ or $\Cal Ri(T)$.
Whether $T$ can be chosen in the class of rank-one actions?
\item
In view of  Theorem~2.1 and Remark~2.2,  whether $\Cal R(\widetilde T)=\Cal E\Cal R(T)$ in the non-Abelian case?
\item
Does  Corollary 2.3 extends to the non-Abelian case, i.e. whether $\Cal E\Cal R(T)=\Cal R(\widehat T)$, where $\widehat T$ is an extension of $T$ to a $G$-action on the same measure space where $T$ is defined? 
\item
A multiple recurrence (and even recurrence) along directions can be defined in the following way.
Let $T$ be a measure preserving action of $\Gamma$ on a $\sigma$-finite measure space $(X,\mu)$ and let $p\in\Bbb N$.
We call $T$ {\it $p$-recurrent along  a line} $\theta\in P(\goth g)$ if for each $\epsilon>0$ and every subset $A\subset X$ of positive measure, there is an element $\gamma\in\Gamma\setminus\{1_\Gamma\}$ and an element $g\in\exp(\theta)$ such that $\text{dist}(\gamma,g)<\epsilon$
and $\mu(A\cap T_\gamma A\cap\cdots \cap T_\gamma^p A)>0$.
Denote by $\Cal R_p(T)$ the set of all $\theta\in P(\goth g)$ such that $T$ is $p$-recurrent along $\theta$. 
Then $\Cal R(T)=\Cal R_1(T)\supset\Cal R_2(T)\supset\cdots$ and $\bigcap_{p\ge 1}\Cal R_p(T)\supset\Cal Ri(T)$.
We note that all these  inclusions are strict
and every set $\Cal R_p(T)$ is  a $G_\delta$.
The results obtained in this work for $\Cal R(T)$ extends to $\Cal R_p(T)$ with similar proofs for each $p$.
\item
Let $T$ be a $(C,F)$-action of $\Gamma$ associated with a sequence $(C_n,F_{n-1})_{n\ge 1}$ satisfying (I)--(IV) and \thetag{3-1} from Section~3.
Given $d>0$, we denote by $C_n^{\otimes d}$ and $F_n^{\otimes d}$ the $d$-th Cartesian power of $C_n$ and $F_n$ respectively.
Then the sequence $(C_n^{\otimes d},F_{n-1}^{\otimes d})_{n\ge 1}$ of subsets in $\Gamma^d$ satisfies (I)--(IV) and \thetag{3-1} from Section~3.
It is easy to see that the $(C,F)$-action $T^{\otimes d}$ of $\Gamma^{ d}$ is canonically isomorphic to the $d$-th tensor product of $T$, i.e. $T^{\otimes d}_{(\gamma_1,\dots,\gamma_d)}=T_{\gamma_1}\times\cdots\times T_{\gamma_d}$ for all $\gamma_1,\dots,\gamma_d\in\Gamma$.
The Lie algebra $\goth g^d$ of $G^d$ is $\goth g\otimes\cdots\otimes\goth g\  (d \text{ times})$.
There is a natural shiftwise action of the permutation group $\Sigma_d$ on $\goth g^d$. This action pushes  down to the protective space $P(\goth g^d)$.
It is easy to see that the sets $\Cal R(T^{\otimes d})$ and $\Cal Ri(T^{\otimes d})$ are invariant under $\Sigma_d$.
In the case where $G=\Bbb R$ and $\Gamma=\Bbb Z$, 
Theorem~4.2 can be refined in the following way: given
a $\Sigma_d$-invariant  subset $\Delta \subset P(\Bbb R^d)$ and a countable $\Sigma_d$-invariant subset  $D$  of $ \Delta$,
 there is a rank-one free infinite measure preserving action $T$ of $\Bbb Z$ such that $D\subset\Cal R(T^{\otimes d})\subset \Delta$.
In particular, each countable $\Sigma_d$-invariant   $G_\delta$-subset  $D$ of $P(\Bbb R^d)$ is realizable as $\Cal R(T^{\otimes d})$ for some  rank-one free action $T$ of $\Bbb Z$.
This generalizes and refines  partly\footnote{This refinement is  partial because we consider only the recurrence set while Adams and Silva studied  simultaneously the set of rational ergodic directions for $T^{\otimes 2}$.} one of the main results from the recent paper by  Adams and Silva \cite{AdSi}: for each $\Sigma_2$-invariant subset $D$ of {\it rational} directions, there is a rank-one action $T$ of $\Bbb Z$ such that $D$ is the intersection of $\Cal R(T^{\otimes 2})$ with the set of all rational directions in $\Bbb R^2$.
We also note that the $\Bbb Z^d$-action $T$  constructed in Theorem~4.3 has the form $T=S^{\otimes d}$ for a $(C,F)$-action $S$ of $\Bbb Z$.
 \item
The theory of directional recurrence can be generalized in a natural way from the infinite measure preserving $\Gamma$-actions to the nonsingular $\Gamma$-actions.
\endroster

\comment

!!!!!!!!!!!!!!!!!!!!!!!!!!!!

\proclaim{Theorem 2.1}
 $\Cal R(\widetilde T)=\Cal E\Cal R(T)$.
\endproclaim
\demo{Proof}
Since $\Gamma$ is uniform,  we may assume
without loss of generality that the subset $\{s(y)\mid y\in G/\Gamma\}$  is relatively compact in $G$.
Hence for each $\epsilon>0$, there is $\epsilon'>0$ such that if $d(g,g')<\epsilon' $
for some $g,g'\in G$
then $d(s(y)gs(y)^{-1},s(y)g's(y)^{-1})<\epsilon$ for each $y\in G/\Gamma$.
Moreover, for each compact $K\subset G$, there is a compact $K'\subset G$ such that  $K\subset\bigcap_{y\in G/\Gamma}s(y)K's(y)^{-1}$.

We first prove  the inclusion $\Cal E\Cal R(T)\subset\Cal R(\widetilde T)$.
Let $\delta\in \Cal E\Cal R(T)$.
Take  an arbitrary subset $B\subset G/\Gamma\times X$ of positive measure.
Then there is a subset $A\subset X$ and an open subset $U\subset G/\Gamma$ such that
$(\lambda\times\mu)(B\cap(U\times A))>(1-\epsilon)(\lambda\times\mu)(U\times A)$
and $s$ is continuous on  $U$.

Since $h_s(1,y)=1$ for all $y\in G/\Gamma$, the standard compactness argument yields that there is an open subset $V\subset U$ and $\vartheta>0$ such that
\roster
\item"(i)"
$gz\in U$
and $h_s(g,z)=1$ for all $z\in V$ and $g\in G$ such that $d(g,1)<\vartheta$
 and
\item"(ii)"
$\lambda(V)>0.9\lambda(U)$.
\endroster
Given a compact $K\subset G$ and a real $\epsilon>0$,
let $K'$ and $\epsilon'$ be determined by $K$ and $\epsilon$ respectively in the aforementioned way.
In a similar way, let $\vartheta'$ be determined by $\vartheta$.
Take  $y_0\in V$.
Since $T$  is evenly recurrent along $\delta$,
it follows that $T$ is recurrent along $s(y_0)^{-1}\cdot\delta$.
Hence
 Lemma~1.5 yields that
 there exist $\gamma_1,\dots,\gamma_n\in \Gamma$, $t_1,\dots,t_n\in\exp(s(y_0)^{-1}\cdot\delta)$ and mutually disjoint subsets $A_1,\dots,A_n$  of $A$ such that $\gamma_j\not\in K'$, $d(\gamma_j,t_j)<\min(\vartheta',\epsilon')$, $T_{\gamma_j} A_j\subset A$ for each $j=1,\dots,n$,
  the sets $T_{\gamma_1} A_1,\dots,T_{\gamma_n} A_n$ are mutually disjoint
  and $\mu(\bigsqcup_{j=1}^nA_j)>0.4\mu(A)$.
We now let $g_{\gamma_j}:=s(y_0)\gamma_j s(y_0)^{-1}$ and
$g_{t_j}:=s(y_0)t_js(y_0)^{-1}$.
Then $g_{t_j}\in\exp(\delta)$, $g_{\gamma_j} y_0=y_0$ and $d(g_{t_j},g_{\gamma_j})<\min(\vartheta,\epsilon)$ for each $j=1,\dots,n$.
Therefore
$$
h_s(g_{\gamma_j},y_0)=s(y_0)^{-1}g_{\gamma_j} s(y_0)=\gamma_j,\qquad j=1,\dots,n.
$$
Hence there is a neighborhood $V_{y_0}$ of $y_0$ such that $V_{y_0}\cup g_{\gamma_j} V_{y_0}\subset V$ and
$h_s(g_{\gamma_j},y)=\gamma_j$ for all $y\in V_{y_0}$, $j=1,\dots,n$.
Since $d(g_{t_j}g_{\gamma_j}^{-1},1)=d(g_{t_j},g_{\gamma_j})<\vartheta$, it follows from~(i) that
\roster
\item"(iii)"
 $g_{t_j}y=g_{t_j}g_{\gamma_j}^{-1}g_{\gamma_j}y\in g_{t_j}g_{\gamma_j}^{-1}V\subset U$
for all
$y\in V_{y_0}$  and
\item"(iv)"
$h_s(g_{t_j}g_{\gamma_j}^{-1},y)=1$ for all $y\in V$.
\endroster
 It follows from (iv) that
$$
h_s(g_{t_j},y)=h_s(g_{t_j}g_{\gamma_j}^{-1}g_{\gamma_j},y)=h_s(g_{t_j}g_{\gamma_j}^{-1},g_{\gamma_j} y)h_s(g_{\gamma_j},y)=h_s(g_{\gamma_j},y)=\gamma_j
$$
for all $y\in V_{y_0}$, $j=1,\dots,n$.
Hence
 $\widetilde T_{g_{t_j}}(y,x)=(g_{t_j}y, T_{\gamma_j} x)$ for all $y\in V_{y_0}$,  $x\in X$ and $j=1,\dots,n$.
This and (iii) yield that
  $$
    \widetilde T_{g_{t_j}} (V_{y_0}\times A_j)
 =(g_{t_j} V_{y_0})\times (T_{\gamma_j} A_j)\subset U\times A,\qquad j=1,\dots,n.
  $$
  The sets $\widetilde T_{g_{t_j}}(V_{y_0}\times A_j)$, $j=1,\dots,n$, are pairwise disjoint.
  Next, $g_{\gamma_j}\not\in K$ and hence $g_{t_j}$ does not belong to the $\epsilon$-neighborhood of $K$.
  Thus we can define a Borel bijection $R_{y_0}$

  It follows that there are finitely many elements $y_1,\dots,y_m\in G/\Gamma$, $g_{t_1},\dots,g_{t_m}\in\delta$ and $\gamma_1,\dots,\gamma_m\in\gamma$ such that $\widetilde T_{g_{t_j}}(V_{y_j}\times (A\cap T_{\gamma_j}^{-1}A))\subset U\times A$ for each $j$ and $\lambda(\bigcup_{j=1}^mV_{y_j})\ge 0.5\lambda (V)$.
  We now let $X_1:=(V_{y_1}\times (A\cap T_{\gamma_1}^{-1}A))$ and
  $X_j:=(V_{y_j}\setminus\bigcup_{i=1}^{j-1}V_{y_j})\times (A\cap T_{\gamma_j}^{-1}A)$ for $j=2,\dots,m$.
We now define  a one-to-one Borel map
  $$
  S:\bigsqcup_{j=1}^mX_j\to \bigsqcup_{j=1}^m \widetilde T_{g_{t_j}}X_j
  $$
  by setting $S\restriction X_j:=\widetilde T_{g_{t_j}}$, $j=1,\dots,m$.

 Hence
$$
\align
(\lambda\times\mu)(\widetilde T_{g_t}B\cap B)
&\ge
(\lambda\times\mu)(\widetilde T_{g_t}(B\cap (U\times A))\cap B\cap(U\times A))\\
&\ge
(\lambda\times\mu)(\widetilde T_{g_t}(U\times A)\cap(U\times A))
-2\epsilon(\lambda\times\mu)(U\times A)\\
&=
\lambda(g_tU\cap U)\mu(T_\gamma A\cap A)-2\epsilon\lambda(U)\mu(A)\\
&>\frac 12\lambda(U)(\mu(T_\gamma A\cap A)-4\epsilon\mu(A))\\
&>0.
\endalign
$$
This means that $\delta\in\Cal R(\widetilde T)$, as desired.

We now show the converse inclusion $\Cal E\Cal R(T)\supset\Cal R(\widetilde T)$.
Let $\delta\in \Cal R(\widetilde T)$.
Take $\epsilon>0$ and a subset $A\subset X$ of positive measure.
Let $U$ be a neighborhood of $\Gamma$ in $G/\Gamma$ such that $s$ is a homeomorphism when restricted to $U$ and
$$
d(s(y),1_G)\le\epsilon\qquad\text{for all $y\in U$.}\tag2-1
$$
Let $V\subset U$ be a neighborhood of $\Gamma$ in $G/\Gamma$ such that
if $hV\cap V\ne\emptyset$ for some $h\in G$ then $h\Gamma\in U$.
Since $\widetilde T$ is recurrent along $\delta$, for each $\epsilon>0$, there is $g\in\exp(\delta)$
such that
$$
(\lambda\times\mu)((V\times A)\cap {\widetilde T}_g(V\times A))>0\tag2-2
$$
 and $g\not\in\bigcup_{y\in G/\Gamma} s(y)$.
Represent $g$ as the product $g=s(g\Gamma)\gamma$ for  some $\gamma\in\Gamma$.
By our choice of $g$, we have $\gamma\ne 1_\Gamma$.
Then
$$
h_s(g,\Gamma)=s(g\Gamma)^{-1}gs(\Gamma)=\gamma.\tag 2-3
$$
It follows from \thetag{2-2} that
 $gV\cap V\ne\emptyset$.
 Hence  $g\Gamma\in U$.
 Since $s$ is continuous on $U$, we obtain that
 the map $V\ni y\mapsto h_s(g,y)\in\Gamma$ is continuous.
Therefore  \thetag{2-3} implies that $h_s(g,y)=\gamma$ for all $y\in V$.
Now \thetag{2-2}  yields that $\mu(T_\gamma A\cap A)>0.$
Since $d$ is right-invariant, it follows from \thetag{2-1} that  $d(g,\gamma)=d(s(g\Gamma),1_G)<\epsilon$.
Hence $\delta\in\Cal R(T)$.
Since by Proposition~1.2,  $\Cal R(\widetilde T)$ is  $G$-invariant, we deduce that
 $\delta\in\Cal E\Cal R(T)$.
\qed
\enddemo

\endcomment

\comment

\remark{Remark \rom{10}}
\roster
 \item
 Let $T$ extend to an action $\widehat T$ of $G$ on the same space.
 Is it true that $\Cal R(T)=\Cal R(\widehat T)$?
 I think no.
 \item
 The set of all actions $T$ of $\Gamma$  such that $\Cal R(T)=P(\goth g)$ is a dense $G_\delta$.
\endroster
We note that (1) follows from (2) and Theorem~9.
\endremark

Let $(Y,\goth Y,\nu,Q)$ be another dynamical system and let $p:X\to Y$ be a nonsingular map that intertwines $T$ and $Q$.
In other words, $Q$ is a factor of $T$ (not necessarily $\sigma$-finite).
Then $\Cal R(T)\subset\Cal R(Q)$.

\proclaim{Proposition 11}
If the extension $T\to Q$ is relatively finite measure preserving then
 $\Cal E\Cal R(T)=\Cal E\Cal R(Q)$.
\endproclaim

\demo{Proof}
Let $\widetilde T$ and $\widetilde Q$ stand for the  $G$-actions induced by $T$ and $Q$ respectively.
Then there is a canonical map intertwining $\widetilde T$ with $\widetilde Q$ \cite{Zi}.
Moreover, the extension $\widetilde T\to\widetilde Q$ is r.f.m.p. whenever the underlying extension $T\to Q$ is so.
It remains to apply Theorem~5, Lemma~7 and the following proposition which is a slight generalization of \cite{Aa, Proposition~1.2.4.}.
\qed
\enddemo

\proclaim{Proposition 12} Let $V$ be a nonsingular transformation and $W$ a r.f.m.p. extension of $W$.
If $V$ is conservative then $W$ is conservative.
\endproclaim

\endcomment

 \comment
 The subset $ \Cal A^G_{\mu,\text{free}}$ of free ($\mu$-almost everywhere) actions is an invariant  $G_\delta$-subspace of  $ \Cal A^G_\mu$.
 Given $A\in \Cal A^G_{\mu,\text{free}}$, the  \rom{Aut}$(X,\mu)$-orbit of $A$ is dense in
 $ \Cal A^G_\mu$.
 \endcomment

\comment

\proclaim{Lemma} Let $\theta\in\Cal R(T)$.
Then for each $\epsilon>0$, $\delta>0$ and a subset $A\subset X$ of positive measure,
there are a subset $A_0\subset A$ and an element $\gamma\in\Gamma$ such that $\mu(A_0)>0$, $d(\gamma,\exp(\theta))<\epsilon$ and $T_\gamma A_0\subset A$ and $1\le \frac{d\mu\circ T_\gamma}{d\mu}(x)<1+\delta$ for all $x\in A_0$.
\endproclaim

\demo
Suppose that the contrary holds.
Then there are $\epsilon>0$, $\delta>0$ and a subset $A$ such that for each $\gamma\in\Gamma$ with $\mu(T_\gamma A\cap A)>0$ and $d(\gamma,\exp(\theta))<\epsilon$, we have
$|\frac{d\mu\circ T_\gamma}{d\mu}(x)-1|>\delta$ for a.a. $x\in T_\gamma A\cap A$.
\enddemo

\endcomment

\comment
 If $T\in \Cal A^\Gamma_{\mu,\text{free}}$ then
$\widetilde T\in  \Cal A^G_{\lambda_{G/\Gamma}\times\mu,\text{free}}$.
\endcomment

\comment
 \proclaim{Lemma 4.5} Let $T\in \Cal A^\Gamma_\mu$.
 If $\theta\not\in\Cal E\Cal R(T)$ then for each $1\ne g\in\exp\theta$, the transformation $(\widetilde T_h)^\bullet$ has positive entropy.
 \endproclaim
 \demo{Proof}
 If $\theta\not\in\Cal E\Cal R(T)$ then $\theta\not\in\Cal R(\widetilde T)$ by Theorem~2.1.
 Then for each $1\ne g\in\exp\theta$, the transformation $\widetilde T_g$ is not conservative by Lemma~2.3(i).
 Let $(\widetilde X,\widetilde \mu)$ stand for the space of $\widetilde T$.
 Then $\widetilde X$ partitions into two $\widetilde T_g$-invariant subsets $X_d$ and $X_c$ of positive measure such that the restriction of $\widetilde T_g$ to $X_d$ is totally dissipative and the restriction of $\widetilde T_g$ to $X_d$ is conservative (see \cite{Aa}).
 Then the Poisson suspension  $\widetilde T_g^\bullet$ is canonically isomorphic to the product $\widetilde (T_g|X_d)^\bullet\times(\widetilde T_g|X_c)^\bullet$ of Poisson suspensions of $T_g|X_d$ and $T_g|X_c$ \cite{Ro}.
 Moreover,  $\widetilde (T_g|X_d)^\bullet$ is Bernoulli \cite{Ro}.
 Hence
 $$
 h( (\widetilde T_g)^\bullet)\ge  h((\widetilde T_g|X_d)^\bullet)>0.
 $$
 \qed
 \enddemo
 \endcomment

\comment

 \proclaim{Theorem 4.6} If there is $T\in  \Cal A^\Gamma_{\mu,\text{free}}$ such that
 $\Cal R(T)=P(\goth g)$ then there is a dense $G_\delta$-subset of $\Cal A^\Gamma_{\mu,\text{free}}$ with this property.
 \endproclaim
 \demo{Proof}
 Let $\Cal Z$ be as in Lemma~4.2.
 It follows from Lemmata 4.2--4.4 that the subset
 $$
\Cal W:= \{T\in  \Cal A^\Gamma_{\mu}\mid (\widetilde T)^\bullet\in\Cal Z\}
 $$
 is an invariant $G_\delta$ in
 $ \Cal A^\Gamma_{\mu}$.
 By Lemma~4.5,
  $\Cal V:=\{T\in \Cal A^\Gamma_{\mu}\mid\Cal R(T)=P(\goth g)\}\supset \Cal W$.
  It now follows from the condition of the lemma and Lemma~4.1 that $\Cal W$ is dense.
  \qed
 \enddemo

 \proclaim{Corollary 4.7}
 There is a dense $G_\delta$ subset of $\Bbb Z^m$-actions which are recurrent along every line.
  \endproclaim

\endcomment

\comment

\example{Example 23}
Let $C_n=\{c(v_nj)\mid 0\le j<N_n\}\cup\{a(v_n)c(v_nj)\mid 0\le j<N_n\}$ if $n$ is odd and
$C_n=\{c(v_nj)\mid 0\le j<N_n\}\cup\{b(v_n)c(v_nj)\mid 0\le j<N_n\}$ if $n$ is even, where $N_n$ and $v_n$ are huge numbers.
Then $\Cal R(T)=\{l_n,l_n'\mid n=0,1,\dots\}$ where $l_n$ is as above and $l_n'$ is a line such that $\exp(l_n')\ni c(n)b(1)$.
Indeed, given odd $n$, $f\in F_m$, $d\in C_{m+1}\cdots C_{n-1}$ and $j<N_n$, we have
$$
T_{a(v_n)}[fdc(v_nj)]_{n}=[fda(v_n)c(v_n(j+j'))]_{n},
$$
where $j'$ is defined by $a(v_n)fda(v_n)^{-1}=fdc(v_nj')$.
Moreover, the union of all $[fdc(v_nj)]_{n}$, $d\in C_{m+1}\cdots C_{n-1}$, $0\le j<N_n$ has a measure which is a half of the measure of $[f]_m$ and $T_{a(v_n)}[fdc(v_nj)]_{n}\subset[f_m]$.
Hence $l_0\in P(\goth h_3(\Bbb R))$.
In a similar way we can verify that $l_n,l_n'\in P(\goth h_3(\Bbb R))$ for all $n$.
\endexample

\endcomment

\comment

Let $G=H_3(\Bbb R)$.
Let $T$ be the $(C,F)$-action of $\Bbb R^2$ associated with a sequence $(C_n,F_{n-1})_{n=1}^\infty$.
Let $F_n$ be a rectangular in $\Bbb R^2$.
We define a sequence of subsets in $H_3(\Bbb R)$ as follows.
Let $\widetilde C_n:=\{b(t_2)a(t_1)\mid (t_1,t_2)\in C_n\}$ and $\widetilde F_n:=\{c(t_3)b(t_2)a(t_1)\mid (t_1,t_2)\in F_n, t_3\in I_n\}$, where
$I_1\subset I_2\subset\cdots$ is a sequence of intervals in $\Bbb R$ such that the sequence $(\widetilde C_n,\widetilde F_{n-1})_{n=1}^\infty$ satisfies the properties of $(C,F)$-actions.
Let $\widetilde T$ be the $(C,F)$-action of $H_3(\Bbb R)$ associated with the sequence $(\widetilde C_n,\widetilde F_{n-1})_{n=1}^\infty$.
Let $X$  be the space of $T$ and let $\widetilde X$ be the space of $\widetilde T$.
Then  a Borel map $\tau:X\to\widetilde X$ is well defined by the following formula:
$$
\tau(f_n,c_{n+1}, c_{n+2},\dots)\mapsto
(p(f_n),p(c_{n+1}), p(c_{n+2}),\dots).
$$
Moreover, $\tau(T_gx)=\widetilde T_{p(g)}\tau(x)$.
It is easy to see that $\tau$ is the $\{\widetilde T_{c(t)}\mid t\in\Bbb R\}$-ergodic decomposition.

 \endcomment

\comment

\example{Simple example} Let $(X,\mu)=(\Bbb R,\lambda)$.
Fix an irrational $\theta$.
Define an action $T$ of $\Bbb Z^2$ by setting $T_{n,m}x=x+n+m\theta$.
Then $\Cal R(T)$ is the line  $\{(t,\theta t)\mid t\in\Bbb R\}$.
\endexample

\subhead Locally compact extensions
\endsubhead
Let $\alpha$ be a cocycle of $T$ with values in a  locally compact second countable group $K$.

\proclaim{Proposition}
Let $S=(S_g)_{g\in G}$ be an ergodic action of $G$ on a standard $\sigma$-finite measure space.
Let $R$ be an invertible transformation commuting with $S_g$ for each $g\in G$.
If $R$ is not conservative then it is totally dissipative.
\endproclaim

\definition{Definition 17}
Given $l\in P(\goth g)$, $\alpha$ is called $l$-recurrent if for each neighborhood $U$ of the identity in $K$, $\epsilon>0$ and a subset $A\subset X$ of positive measure, there is $\gamma\in \Gamma$
and a  subset $A_0\subset A$ such that $\gamma\ne\text{Id}$, $\mu(A_0)>0$, $\gamma A_0\subset A$, $d(l,\gamma)<\epsilon$  and $\alpha(x,\gamma)\in U$ for each $x\in A_0$.
\enddefinition

\proclaim{Proposition 18}
$\alpha$ is $l$-recurrent if and only if the skew product $T^\alpha$ is $l$-recurrent.
\endproclaim

Denote by $\Cal R(\alpha)$ the set of all $l$ such that $\alpha$ is $l$-recurrent.
By $\Cal E\Cal R(\alpha)$ we denote the set of all $l$ such that the $G$-orbit of $l$ is in $\Cal R(\alpha)$.
Of course, $\Cal R(\alpha)$ and $\Cal E\Cal R(\alpha)$ are cohomology invariants, $\Cal R(\alpha)\subset\Cal R(T)$,
$\Cal E\Cal R(\alpha)\subset\Cal E\Cal R(T)$.
The next claim follows from Proposition~11.

\proclaim{Proposition 19} If $K$ is compact then
$\Cal E\Cal R(\alpha)=\Cal E\Cal R(T)$.
\endproclaim

\endcomment

\comment
\demo{Proof}
Let $l\in\Cal R(T)$.
Take a subset $A$ of positive measure, a neighborhood $U$ of the identity in $K$ and $\epsilon>0$.
Select a symmetric neighborhood $V$ go the identity in $G$ such that $VVV\subset U$.
Then we construct a sequence of subsets $A_1\supset A_2\supset\cdots$
a sequence of elements $\gamma_1,\gamma_2,\dots$ in $\Gamma$, a sequence $k_1,k_2,\dots$ of elements in $K$, a sequence of positive reals $\epsilon_1,\epsilon_2, \dots$ such that 	for each $j>0$, we have
\roster
\item
 $A_j\cup T_{\gamma_j} A_j\subset T_{\gamma_{j-1}}A_{j-1}$,
 \item
  $d(l,\gamma_j)<\epsilon_j$ and
  \item
  $\alpha(x,T_{\gamma_j})\in Vk_j$ for all $x\in A_j$.
  \endroster
  We assume that $A_0:=A$ and $\gamma_0:=1_\Gamma$.
  Since $K$ is compact, there are integers $j_1<j_2$ such that $k_{j_2}k_{j_1}^{-1}\in V$.
  We have
  $$
  T_{\gamma_{j_2}}A_{j_2}\subset T_{\gamma_{j_2-1}}A_{j_2-1}\subset\cdots\subset
  T_{\gamma_{j_1}}A_{j_1}
  $$
  and hence $T_{\gamma_{j_1}^{-1}\gamma_{j_2}}A_{j_2}\subset A_{j_1}\subset A$.
  We see that for each $x\in A_{j_2}$,
  $$
  \align
  \alpha(x,T_{\gamma_{j_1}^{-1}\gamma_{j_2}})&=\alpha(x, T_{\gamma_{j_2}})\alpha(T_{\gamma_{j_2}}x, T_{\gamma_{j_1}^{-1}})\\
 & =
  \alpha(x, T_{\gamma_{j_2}})\alpha(T_{\gamma_{j_1}^{-1}}T_{\gamma_{j_2}}x, T_{\gamma_{j_1}})^{-1}
  \\
 & \in Vk_{j_2}k_{j_1}^{-1}V\subset VVV\subset U.
  \endalign
  $$
\enddemo
\endcomment

\comment

 \proclaim{Corollary 4.9}
 There is a Poisson action $V$ of $H_3(\Bbb R)$   of $0$ entropy such that for each $1\ne g\in H_3(\Bbb R)$, the entropy of  $V_g$ is strictly positive.
\endproclaim
 \demo{Proof}
 Let $T$ be as in Proposition~3.6.
 Since $T$ is of rank one, the entropy of $T^\bullet$ is trivial.
 The Poisson  suspension of $\widetilde T$ is the co-induced action of $T^\bullet$.
The entropy of the co-induced action is trivial whenever the entropy of $T^\bullet$ is trivial.
\qed
 \enddemo

\endcomment


\Refs
 \widestnumber\key{Au--Ru}


 \ref
 \key Aa
 \by J. Aaronson
 \book
 Introduction to infinite ergodic theory
 \yr 1997
\vol 50
\publ Amer. Math. Soc.
\publaddr
Providence, R.I.
\bookinfo  Mathematical Surveys and Monographs
 \endref


\ref
\key AdSi
\by T.  M. Adams and C. E. Silva
\paper
On infinite transformations with maximal control of ergodic two-fold product powers
\paperinfo
ArXiv: 1402.1818v1
\endref



\ref
\key Au--Ha
\by L. Auslander, L. W. Green and F. Hahn
\book  Flows on homogeneous spaces
\publ
 Princeton University Press
 \publaddr
  Princeton, N.J.
  \yr 1963
\endref






 \ref
 \key Da1
 \by A. I. Danilenko
 \paper
$(C, F)$-actions in ergodic theory
\inbook
 Progr. Math.
 \vol  265
 \publ
 Birkh\"auser
 \publaddr
  Basel
 \yr 2008
 \pages 325--351
\endref


\ref
\key Da2
\bysame
\jour  Erg. Th. \& Dynam. Sys.
\paper
Mixing actions of the Heisenberg group
\vol
 34
 \yr 2014
 \pages 
  1142--1167
\endref



 \ref
 \key Da3
 \paperinfo arXiv: 1403.2012
 \paper    Actions of finite rank: weak rational ergodicity and partial rigidity
   \bysame
   \endref


\ref
\key DaLe
\by A. I. Danilenko and M. Lema{\'n}czyk
\paper
Odometer actions of the Heisenberg group
\jour J. d'Anal. Math.
\toappear
\endref

\ref
\key DaSi
\inbook Encyclopedia of Complexity and Systems Science
\yr 2009
\pages 3055--3083
\by    A. I. Danilenko and 
    C. E. Silva 
\paper  Ergodic theory:
nonsingular transformations
\publ
Springer
\publaddr
New York
\endref


\ref\key Fel
\by J. Feldman
\paper A ratio ergodic theorem for commuting, conservative, invertible transformations with quasi-invariant measure summed over symmetric hypercubes
\jour Erg.
Th. \& Dynam. Sys.
\vol 27
\yr 2007
\pages 1135--1142
\endref

 \ref
 \key FeKa
 \by    S. Ferenczi and B. Kami{\'n}ski
    \paper Zero entropy and directional Bernoullicity of a Gaussian $\Bbb Z^{2}$-action
\jour
    Proc. Amer. Math. Soc.
    \vol 123
    \yr 1995
    \pages 3079--3083
    \endref


 \ref
\key  Ja--Ru
\by E. Janvresse, T. Meyerovitch, E. Roy and T. de la Rue
\jour Trans. Amer. Math. Soc.
\paper
Poisson suspensions and entropy for infinite transformations
\vol 362
\yr 2010
\pages  3069--3094
\endref


 \ref
 \key JoSa
 \paper Directional recurrence for infinite measure preserving $\Bbb Z^d$-actions
 \jour Erg. Th. \& Dynam. Sys.
 \toappear
\by A. S. A. Johnson and A. A. {\c S}ahin
  \endref

 \ref
 \key Ma
 \by
 G. Mackey
 \paper
 Induced representations of locally compact groups. I
 \jour Ann. Math.
 \vol 55
 \yr 1952
 \pages 101--139
 \endref



 \ref
 \key Mal
 \by A. I. Malcev
 \paper On a class of homogeneous spaces
 \jour  Izvestiya Akad. Nauk. SSSR. Ser. Mat.
 \vol 13
 \yr 1949
 \pages 9--32
 \lang Russian
 \endref

 \ref
 \key Mi
 \by J. Milnor
\paper On the entropy geometry of cellular automata
\jour Complex Systems
\vol 2
\yr 1988
\pages 357--385
\endref


\ref
\key Ne
\by Yu. Neretin
\book
 Categories of symmetries and infinite-dimensional groups
  \bookinfo
  London Mathematical Society Monographs. New Series, 16. Oxford Science Publications. The Clarendon Press, Oxford University Press, New York, 1996
  \endref
  
  
  \ref\key Pa
  \by K. K. Park
  \jour  Isr. J. Math. 
  \vol 113 
  \yr 1999
  \pages  243--267
  \paper On directional entropy functions
  \endref


 \ref
 \key Ro
 \by E. Roy
 \paper
 Poisson suspensions and infinite ergodic theory
 \jour Erg. Th. \& Dynam. Sys.
\vol 29
\yr 2009
\pages 667--683
 \endref

 \ref
 \key Zi
 \by R. J. Zimmer
 \paper Induced and amenable ergodic actions of Lie
groups
 \jour Ann. Sci. Ecole Norm. Sup.
\vol 11
\yr 1978
\pages 407--428
\endref
\endRefs
\enddocument